\documentclass[11pt]{article}
\usepackage{amsfonts}
\usepackage{amsmath,amsthm,amscd,amssymb,mathrsfs,setspace, textcomp,bm}
\usepackage{enumerate}
\usepackage{latexsym,epsf,epsfig}
\usepackage{color}
\usepackage[hmargin=2.25cm,vmargin=2.75cm]{geometry}

\usepackage{hyperref}
\usepackage{cite}

\newcommand{\ds}{\displaystyle}
\newcommand{\db}{\Big|\Big|}
\newcommand{\reals}{\mathbb{R}}
\newcommand{\N}{\mathbb{N}}
\newcommand{\realstwo}{\mathbb{R}^2}
\newcommand{\realsthree}{\mathbb{R}^3}
\newcommand{\nn}{\nonumber}
\newcommand{\dt}{{[\Delta t]}}
\newcommand{\xb}{{\bf{x}}}
\newcommand{\bfn}{{\bf{n}}}
\newcommand{\bfeta}{{\boldsymbol{\eta}}}
\newcommand{\bfpsi}{{\boldsymbol{\psi}}}
\newcommand{\bfsigma}{{\boldsymbol{\sigma}}}
\newcommand{\Dt}{{\partial_t}}
\newcommand{\Dx}{{\partial_x}}
\newcommand{\cA}{{\mathcal{A}}}
\newcommand{\cH}{\mathcal{H}}
\newcommand{\ccH}{\mathscr{H}}
\newcommand{\cD}{\mathscr{D}}
\newcommand{\Dn}{\partial_{\nu}}
\newcommand{\DT}{\partial_{\tau}}
\newcommand{\Dz}{\partial_z}
\newcommand{\cE}{{\mathcal{E}}}
\newcommand{\bA}{\mathbf A}
\newcommand{\bV}{\mathbf V}
\newcommand{\bF}{\mathbf F}
\newcommand{\cF}{\mathcal{F}}
\newcommand{\cG}{\mathcal{G}(z)}
\newcommand{\cP}{\mathcal{P}}
\newcommand{\R}{\mathbb{R}}
\newcommand{\cB}{\mathcal{B}}
\newcommand{\bu}{\mathbf u}
\newcommand{\bw}{\mathbf w}
\newcommand{\bv}{\mathbf v}
\newcommand{\bt}{\boldsymbol t}
\newcommand{\bz}{\mathbf z}
\newcommand{\bp}{\mathbf{\Phi}}
\newcommand{\bE}{\mathbf{E}}
\newcommand{\Om}{\Omega}
\newcommand{\om}{\omega}
\newcommand{\ga}{\gamma}
\newcommand{\s}{\sigma}
\newcommand{\e}{\epsilon}
\newcommand{\Ez}{E_{z}}
\newcommand{\bT}{\mathbb{T}}
\newcommand{\lb}{ \langle}
\newcommand{\rb}{ \rangle}
\newcommand{\ou}{\overline u}
\newcommand{\Lo}{\color{magenta}}
\theoremstyle{plain}
\newtheorem{theorem}{Theorem}[section]
\newtheorem{lemma}[theorem]{Lemma}
\newtheorem{proposition}[theorem]{Proposition}

\newtheorem{definition}{Definition}

\newtheorem{corollary}[theorem]{Corollary}
\newtheorem{assumption}{Assumption}
\theoremstyle{remark}
\newtheorem{remark}{Remark}[section]
\numberwithin{equation}{section}
\numberwithin{theorem}{section}
\numberwithin{remark}{section}
\numberwithin{assumption}{section}
\numberwithin{condition}{section}
\date{}

\begin{document}

\title{Multilayered Poroelasticity Interacting with Stokes Flow}
\author{Lorena Bociu\footnote{North Carolina State University, \url{lvbociu@ncsu.edu}};~~~~Sun\v{c}ica \v{C}ani\'{c}\footnote{University of California-Berkeley, \url{canics@berkeley.edu}};~~~~Boris Muha\footnote{University of Zagreb, Croatia, \url{borism@math.hr}};~~~~Justin T. Webster\footnote{University of Maryland, Baltimore County, Baltimore, MD, \url{websterj@umbc.edu}}}
     \maketitle
         \abstract{
         
        We consider the interaction between an incompressible, viscous fluid modeled by the dynamic Stokes equation and a multilayered poroelastic structure which consists of a thin, linear, poroelastic plate layer (in direct contact with the free Stokes flow) and a thick Biot layer. The fluid flow and the elastodynamics of the multilayered poroelastic structure are fully coupled across a fixed interface
 through physical coupling conditions (including the Beavers-Joseph-Saffman condition), which present mathematical challenges related to the regularity of associated velocity traces.
         We prove existence of weak solutions to this fluid-structure interaction problem with either (i) a linear, dynamic Biot model, or (ii) a nonlinear quasi-static Biot component, where the permeability is a nonlinear function of the fluid content (as motivated by biological applications). The proof is based on constructing approximate solutions
through Rothe's method, and using energy methods and a version of Aubin-Lions compactness lemma (in the nonlinear case) to recover the weak solution as the limit of approximate subsequences. We also provide uniqueness criteria and show that constructed weak solutions are indeed strong solutions to the coupled problem if one assumes additional regularity.

\if 1 = 0         
           We address existence and uniqueness of weak solutions for Stokes flow interacting with a thick poroelastic region, separated by a thin poroelastic plate (utilizing the recent model \cite{mikelic}). The dynamics of poroelastic filtration---the coupling of a Biot-type poroelastic solid coupled to a fluid flow---have been studied before, both theoretically \cite{Sho05,Ce17} and numerically \cite{badia2009coupling,AKY18,AEN18}. This type of coupled model presents challenges in the dynamics at the interface, particularly, the regularity of associated velocity traces. In this treatment, we also consider a slip-dissipative coupling condition known as the {\em Beavers-Joseph-Saffman}  \cite{...}  condition between the Stokes flow and the poroelastic plate. The fully dynamic model (allowing for inertia and compressibility) is considered in the fully linear case (as in \cite{Sho05}), while biologically motivated \cite{bgsw, MBE2,BW} nonlinearity in the permeability tensor is considered for the quasi-static case where inertia may be neglected in the thick Biot region. We show equivalence of the strong form of the coupled dynamics and our particular, non-obvious weak form. Requiring only the constrained storage coefficient of the plate to be strictly positive, weak solutions are defined and explicitly constructed through semi-discretization in time (Rothe's method). We provide uniqueness criterion in the linear case, and a weak-strong uniqueness type result in the nonlinear, quasi-static case.
           \fi
           }

\section{Introduction and Literature Review}
Poroelasticity refers to fluid flow, as described by Darcy's law, within deformable, porous media. Historically, studies of the fluid flow through porous media were motivated by applications in geosciences, see e.g. \cite{detournay1993fundamentals}.  The classic partial differential equation (PDE) model for such flows, based on consolidation theory, was originally introduced by Biot \cite{biot,Biot55,Biot63} and predicated upon phenomenological modeling. 
An alternative approach considers the fluid-structure interaction system at  the level of the pores, 
with a derivation of an effective model based on homogenization techniques, see \cite{Auriault97,FM03,MikWhe12} and the references therein. 
Derivation of poroelastic models based on a continuum mechanics approach can be found in the monograph of Coussy \cite{coussy}.  
Recently, Mikeli\'{c} et al. presented a mathematically rigorous derivation of lower dimensional poroelastic models, 
such as poroelastic plates and shells, which can be found in  \cite{CzoMik15,MikTam16}.
From the point of view of PDE well-posedness, the classical Biot system has been extensively studied by many authors, see, e.g., \cite{zenisek, show1, show2,owc, showmono,barucq1, barucq2}. Since poroelastic structures are ubiquitous in biological systems, recently there has been a growing interest in the analysis of poroelasticity in the context of biomedical applications, see, e.g. \cite{bgsw,cowin1999bone,causin2014poroelastic, huyghe1992porous,MBE, MBE2, BN19, banks1, banks2,guidoboni2006} and the references and discussions therein. 

In a variety of biomedical applications, the poroelastic structure is often in contact with a free fluid flow,
such as, e.g., the blood flow through an artery with poroelastic walls.
Mathematically, such problems are described mathematically by a coupled system of partial differential equations:
the fluid equations (e.g. the Navier-Stokes or Stokes system) on one side, and the equations of poroelasticity on the other. These coupled problems are referred to as fluid-poroelastic structure interaction (FPSI) problems, and have been analyzed, for instance, in \cite{Sho05, Ce17, AEN18,rohan1,rohan2}.   
However, many biological tissues (such as arterial walls) have a {\em multilayered structure}. 
{{In bioengineering, a multilayered poroelastic structure arises in the design of a bioartificial pancreas \cite{DesaiPancreasReview}.
In bioartificial pancreas, healthy pancreatic cells are seeded in a poroelastic agarose or alginate gel, 
and the gel is encapsulated between two poroelastic plates, which are designed to exclude the patients own immune cells from
attacking the organ while allowing passage of oxygen and nutrients to the pancreatic cells \cite{DesaiPancreasReview}. 
This multilayered poroelastic medium is connected to the blood flow via a tube (an anastomosis graft). 
The blood flows over the poroelastic plate and it passes into the poroelastic medium containing the transplanted cells,
bringing oxygen and nutrient supply to the cells. 
Thus, the poroelastic plate in this design separates the poroelastic region containing the pancreatic cells on one side,
from the region containing free blood flow on the other.
One of the problems in this bioartificial pancreas design is to understand fluid (blood) flow 
through the poroelastic plate and poroelastic gel containing the cells, and simulate oxygen 
and nutrients supply to the cells.
}}

{{{\bf{Motivated by these types of applications, in the present paper we are interested in the mathematical aspects of the interaction between the flow of an incompressible, viscous fluid and a multilayered poroelastic structure.}} 
The multilayered poroelastic structure under consideration consists of two layers:
a thin poroelastic layer in direct contact with the free fluid flow, and a thick poroelastic  layer sitting  atop the thin layer (see Fig.~\ref{domain}).}} The fluid flow is described by the time-dependent Stokes equations,
while the thin and thick poroelastic layers are described by a poroelastic plate model, and a Biot model, respectively.
The three different physical models (the fluid model, the thin poroelastic structure model, and the thick Biot system) are fully coupled across a fixed interface by physically-motivated interface conditions. 
We mention that interface conditions between porous media (Biot model(s)) and adjacent free fluid flow were rigorously analyzed in \cite{JagMik96,JagMik00}. 
Although slightly different from  \cite{JagMik96}, our interface conditions involve a jump in the pressure between the free fluid and the Biot pressure at the interface.
However, in contrast with the pressure jump derived in \cite{JagMik96} through boundary layer considerations,
the pressure jump in the present paper is a consequence of the presence of a thin poroelastic plate, which serves as the fluid-structure interface.

{\bf In the multilayered poroelastic structure, we consider two distinct scenarios}: (i) a linear, dynamic Biot model, or (ii) a nonlinear quasi-static Biot component where the permeability is a nonlinear function of the fluid content. The second case is motivated by biological applications \cite{Hsu,causin2014poroelastic,detournay1993fundamentals,DESpaper}. 
In particular, our main motivation comes from the recent work \cite{DESpaper}
on vascular prostheses called drug-eluting stents.
Based on mathematical and computational results, together with experiments cited in  \cite{DESpaper},
it was found that the implantation of a vascular stent within the arterial wall, induces a change 
 in vascular tissue permeability, which is affected by the change in porosity
near the end points (edges) of vascular stents. As a consequence, this can influence drug advection, reaction and diffusion
into the vascular tissue, and lower the efficacy of anti-inflammatory drugs in lowering restenosis (or re-closure) of 
coronary arteries near the edges of a stent.
Thus, in the present paper we accommodate the case when {\emph{permeability is not constant, but depends on porosity}}, and, in turn, on the so called {\em Biot fluid content}. We note that the nonlinear permeability (dependent on fluid content) introduces {\emph{significant complications in the analysis}} of the system. Such a nonlinear  model was first considered---from a mathematical point---in \cite{cao}, and later in \cite{bgsw}. 
The  reference \cite{cao} focuses  on the {\it compressible Biot model}, with a permeability that depends on the Biot dilation, and constructs weak solutions through a full spatio-temporal discretization  in the mathematically simplified framework of homogeneous boundary conditions for both fluid pressure and solid displacement. 
The later reference \cite{bgsw} focuses on {\it Biot models with incompressible constituents} and constructs weak solutions following the approach of \cite{zenisek} for both poroelastic and poro-visco-elastic systems. Another main contribution of \cite{bgsw}  is the treatment of non-homogeneous, mixed boundary conditions that are  physically relevant to ophthalmological applications. 

The mathematical literature dealing with FPSI problems is scarce. 
We mention the semigroup approach based on the variational formulation of a linear {\em compressible} Stokes system coupled to a Biot dynamics developed in \cite{Sho05}. The incompressible Stokes case, considered in the present manuscript,
does not seem to be covered by this analysis.
Later, in \cite{Ce17}, the problem of linear coupling (across a fixed interface) 
of the incompressible Navier-Stokes equations and {\em linear} Biot equations was considered, 
and the existence of {\em strong solutions} was proved under a small data assumption.
And finally, we mention a well-posedness result for the stationary non-Newtonian fluid flow coupled with the Biot equation via fixed interface, reported in \cite{AEN18}. 
To the best of our knowledge there are no results related to the well-possednes of the FPSI problems in cases where the poroelastic model is nonlinear, and additionally, consists of multiple layers. 

{\bf In the present paper} we prove the {\emph{existence}} of a weak solution to a fluid-structure interaction problem between the flow of an 
incompressible, viscous fluid, and a multilayered poroelastic structure, where the Biot model in the thick poroelastic layer
may be nonlinear through the dependence of permeability on the porosity. 
We consider two cases: the linear case in which the permeability in the Biot model may be a given function of space and time, 
and the nonlinear, {\em quasi-static} case
in which the permeability in the  Biot model may generally depend on porosity. 
The proof is based on the construction of approximate solutions via Rothe's method, introduced in \cite{ARMA} to study
fluid-structure interaction problems between the flows of incompressible, viscous fluids and {\emph{elastic}} structures. 
The method is extended in the present paper to the poroelastic case, where the thin poroelastic plate serves as a 
fluid-structure interface, and plays the role of a regularizing mechanism in the existence proof, similar to the results in \cite{BorisSunNWE, reg1,reg2}. 
The proof is based on a construction of approximate solutions by semi-discretizing the coupled problem in time,
and ``solving'' a sequence of elliptic problems for each fixed time step $\Delta t$. 
Using energy methods, we show that as $\Delta t \to 0$, approximate solutions converge to a weak solution
of the underlying coupled problem. In the nonlinear case, a compactness argument, based on a version of the 
Aubin-Lions lemma, is used to show strong convergence of subsequences to a weak solution. Finally, for the linear problem we demonstrate uniqueness in a class of solutions with slightly improved regularity than weak solutions. For the nonlinear problem, we show a weak-strong type uniqueness result for the same class.
This weak-strong uniqueness result for the FPSI with the nonlinear, quasi-static Biot model considered here is new,
since even for the nonlinear Biot model alone, a full uniqueness result seems to be lacking. 

We finish this introduction by mentioning that the development of  numerical methods for FPSI problems has been 
a very active research area in the last decade, see e.g. 
\cite{AKY18,badia2009coupling,bukavc2012fluid} and the references therein. 
Many problems remain open in this field, especially related to the design of partitioned schemes for FPSI problems.
In particular, one of the approaches currently under investigation
is the design of a partitioned scheme 
that would be based on the main steps in the constructive proof presented in this manuscript. 

{{We conclude this section by noting that extensions of this work to the moving interface case, in which the fluid domain is not fixed
but moving with the interface itself, requires careful considerations. Such an analysis would be similar, but more complicated than those presented by the 
authors in \cite{ARMA,BorisSunNWE,BCLT}, where FSI with multi-layered structures and elasticity were considered, but without poroelasticity. This is one of the goals of our future research. }}

\if 1 = 0
\documentclass[11pt]{article}
\usepackage{amsfonts}
\usepackage{amsmath,amsthm,amscd,amssymb,mathrsfs,setspace, textcomp,bm}
\usepackage{enumerate}
\usepackage{latexsym,epsf,epsfig}
\usepackage{color}
\usepackage[hmargin=2.25cm,vmargin=2.75cm]{geometry}

\frenchspacing
\usepackage{hyperref}


\setcounter{MaxMatrixCols}{10}

\newcommand{\ds}{\displaystyle}
\newcommand{\db}{\Big|\Big|}
\newcommand{\reals}{\mathbb{R}}
\newcommand{\realstwo}{\mathbb{R}^2}
\newcommand{\realsthree}{\mathbb{R}^3}
\newcommand{\nn}{\nonumber}
\newcommand{\dt}{{[\Delta t]}}
\newcommand{\xb}{{\bf{x}}}
\newcommand{\bfn}{{\bf{n}}}
\newcommand{\bfeta}{{\boldsymbol{\eta}}}
\newcommand{\bfpsi}{{\boldsymbol{\psi}}}
\newcommand{\bfsigma}{{\boldsymbol{\sigma}}}
\newcommand{\Dt}{{\partial_t}}
\newcommand{\Dx}{{\partial_x}}
\newcommand{\cA}{{\mathcal{A}}}
\newcommand{\cH}{\mathcal{H}}
\newcommand{\ccH}{\mathscr{H}}
\newcommand{\cD}{\mathscr{D}}
\newcommand{\Dn}{\partial_{\nu}}
\newcommand{\DT}{\partial_{\tau}}
\newcommand{\Dz}{\partial_z}
\newcommand{\cE}{{\mathcal{E}}}
\newcommand{\bA}{\mathbf A}
\newcommand{\bV}{\mathbf V}
\newcommand{\bF}{\mathbf F}
\newcommand{\cF}{\mathcal{F}}
\newcommand{\cG}{\mathcal{G}(z)}
\newcommand{\cP}{\mathcal{P}}
\newcommand{\R}{\mathbb{R}}
\newcommand{\cB}{\mathcal{B}}
\newcommand{\bu}{\mathbf u}
\newcommand{\bw}{\mathbf w}
\newcommand{\bv}{\mathbf v}
\newcommand{\bt}{\boldsymbol t}
\newcommand{\bz}{\mathbf z}
\newcommand{\bp}{\mathbf{\Phi}}
\newcommand{\bE}{\mathbf{E}}
\newcommand{\Om}{\Omega}
\newcommand{\om}{\omega}
\newcommand{\ga}{\gamma}
\newcommand{\s}{\sigma}
\newcommand{\e}{\epsilon}
\newcommand{\Ez}{E_{z}}
\newcommand{\bT}{\mathbb{T}}
\newcommand{\lb}{ \langle}
\newcommand{\rb}{ \rangle}
\newcommand{\ou}{\overline u}
\newcommand{\Lo}{\color{magenta}}
\theoremstyle{plain}
\newtheorem{theorem}{Theorem}[section]
\newtheorem{lemma}[theorem]{Lemma}
\newtheorem{proposition}[theorem]{Proposition}
\newtheorem{condition}{Condition}
\newtheorem{definition}{Definition}
\newtheorem{conjecture}{Conjecture}
\newtheorem{corollary}[theorem]{Corollary}
\newtheorem{assumption}{Assumption}
\theoremstyle{remark}
\newtheorem{remark}{Remark}[section]
\numberwithin{equation}{section}
\numberwithin{theorem}{section}
\numberwithin{remark}{section}
\numberwithin{assumption}{section}
\numberwithin{condition}{section}
\linespread{1.05}

\begin{document}

\fi


\section{Mathematical Formulation}
\subsection{Spatial Domain}
\begin{figure}[h!]
\begin{center}
           \includegraphics[width = 0.450 \textwidth]{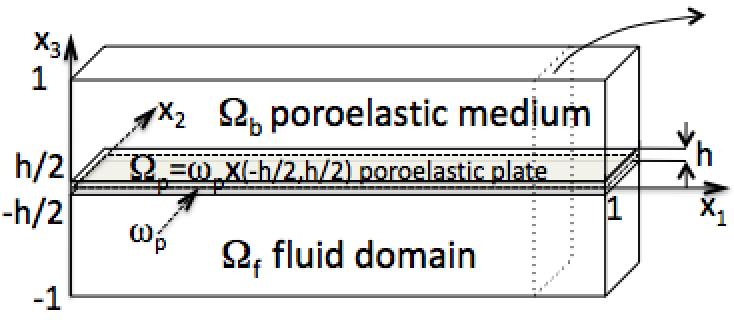}\hskip 0.3in
           \includegraphics[width = 0.260 \textwidth]{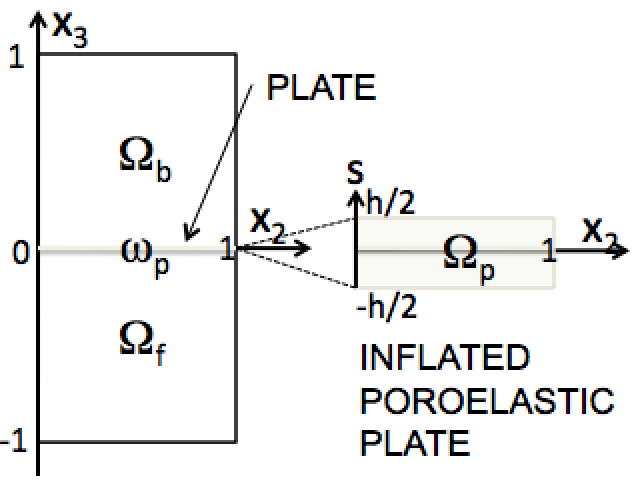}
           \end{center}
        \caption{\small Left: 3D domain $\Omega$. Right: 2D vertical cross-section through domain $\Omega$.}
\label{domain}
\end{figure}

We consider the domain $\Omega\subset \mathbb R^3$ with the simplest geometry describing
the coupling between the flow of an incompressible, viscous fluid and a thick poroelastic medium,
where the coupling occurs through a thin poroelastic plate, serving as a fluid-structure interface with mass.
Namely, our domain is a rectangular prism: $\Omega=(0,1)\times (0,1) \times (-1,1)$, 
with three distinct subregions, see Fig.~\ref{domain}:
$\Omega_f$ which is occupied by a viscous, incompressible fluid, 
$\Omega_b$ which contains a (thick) poroelastic medium modeled by the Biot's model,
and $\Omega_p$ corresponding to a  thin poroelastic plate of thickness $h$, 
whose middle surface, denoted by $\omega_p$, separates the regions $\Omega_f$ and $\Omega_b$: 
\begin{align*}
\Omega_f=&~(0,1)^2\times (-1,0),\ &\Omega_b&=~(0,1)^2\times (0,1), \\
\omega_p=&~(0,1)^2\times \{0\}, \ &\Omega_p&=~(0,1)^2\times  (-h/2,h/2).
\end{align*}
In Cartesian coordinates $(x_1,x_2,x_3)$, the physical domain $\Omega$ is the 
union $\Omega = \Omega_f \cup \omega_p \cup \Omega_b$. 
The variable $x_3$ will be referred to as the ``transverse'' variable $x_3 \in [-1,1]$,
and $x_1\in[0,1]$ and $x_2\in[0,1]$ will be called the ``longitudinal'' and ``lateral'' variables, respectively. 

To study filtration through the thin, poroelastic plate of thickness $h$, 
we will consider a ``$2.5$-dimensional model'', derived and rigorously justified by Marciniak-Caochra and Mikeli\'{c} in \cite{mikelic},
obtained via asymptotic reduction from the full Biot model. 
In the poroelastic plate model, the elastodynamics of the thin plate is described in terms of the $x_1$ and $x_2$ variables,
as usual, 
namely, it is given in terms of the plate displacement of the middle surface,
which is a function of only $x_1$ and $x_2$.
Filtration, however, needs to be described in terms of an additional variable $s\in [-h/2,h/2]$, 
since it  was shown in \cite{mikelic} that filtration through a thin poroelastic isotropic structure is dominant in the thin/transverse 
direction. Therefore, we denoted by $\omega_p$ the location of the middle surface of the elastic plate,
and by $\Omega_p$ the ``inflated'' mathematical domain on which
the filtration velocity and pore fluid pressure will be defined.
In Sec.~\ref{model} we present the model, and provide more details about the approach.


The fluid contained in region $\Omega_f$ will be modeled by the time-dependent Stokes equations, with periodic boundary conditions at 
$x_1 = 0, x_1 = 1, x_2 = 0$, and $x_2 = 1$, and with the no-slip boundary condition at the bottom boundary $x_3 = -1$,
which is considered fixed.
The boundary conditions at the top fluid domain boundary, where the bulk fluid flow in $\Omega_f$  is coupled to the 
filtration flow in through the poroelastic membrane $\Omega_p$, will be described in Sec.~\ref{model}. 

 Finally, $\Omega_b$ represents a thick poroelastic region relative to the thin poroelastic plate, with filtration flow and structure
 displacement satisfying Biot's equations. 
 The thick structure displacement in $\Omega_b$, as well as the
 flitrate flow, will satisfy periodic boundary conditions at $x_1 = 0, x_1 = 1, x_2 = 0$, and $x_2 = 1$. 
 At the top boundary of $\Omega_b$, with the reference configuration described by $x_3 = 1$, 
 we prescribe  zero elastic stress, and no filtration  through the wall.
 At the bottom boundary, both the elastodynamics of the thick poroelastic medium, as well as filtration flow in $\Omega_b$
 are fully coupled to the elastodynamics and filtration flow in the thin poroelastic plate $\Omega_p$,
 which is, in turn, fully coupled to the fluid flow in $\Omega_f$. 
 The filtration flow in $\Omega_b$ is driven by the pressure gradient between
 the bottom boundary $x_3 = 0$, and the top boundary $x_3 = 1$, and by the 
 time-dependent change in structure displacement. 

 In summary, the boundary of $\Omega$ consists of the periodic boundary
 at $x_1 = 0, x_1 = 1, x_2 = 0$, and $x_2 = 1$,
a fixed bottom boundary $x_3 = -1$, and the top boundary with the reference configuration 
described by $x_3 = 1$.
Detailed boundary conditions, and the coupling between different subdomains of $\Omega$,
will be described in the next section.

\subsection{Mathematical Models}\label{model}

In this section we describe the mathematical models, holding in each of the sub-domains specified above,
and specify the coupling between them.
We will be using $\mathbf{x}$ to denote points in Cartesian coordinates $\mathbf{x}=(x_1,x_2,x_3)$,
the subscript $b$ will be used to denote quantities in the Biot domain $\Omega_b$, and subscript $p$
will be used to denote the quantities corresponding to the poroelastic plate.

\subsubsection{Biot Model in $\Omega_b$}
Domain $\Omega_b$ denotes a poroelastic medium, modeled by the following Biot model, given in terms of
the unknown functions
$\bfeta=\bfeta(\mathbf{x},t)$, describing displacement of the solid matrix from its reference configuration $\Omega_b$,
and 
$p_b=p_b(\mathbf{x},t)$ denoting the pore fluid pressure:
\begin{equation}\label{Biot_main}
\begin{cases}
\rho_b\bfeta_{tt}-\mu_b\Delta \bfeta-(\mu_b+\lambda_b)\nabla (\nabla \cdot \bfeta)+\alpha_b \nabla p_b=\mathbf F_b, &\hskip1cm  
{\rm in} \ \Omega_b \times (0,T), \\[.2cm]
[c_bp_b+\alpha_b\nabla \cdot \bfeta]_t-\nabla \cdot (k_b\nabla p_b)=S, &\hskip1cm   {\rm in} \ \Omega_b \times (0,T). \\[.2cm]
\end{cases}
\end{equation}
On the right hand-side, $\mathbf{F}_b=\mathbf{F}_b(\mathbf{x},t)$ denotes a body force per unit volume, and
$S=S(\mathbf{x},t)$ is a net volumetric fluid production rate. 
The first equation describes the balance of linear momentum for the fluid-solid mixture,
and the second  equation describes conservation of mass for the fluid component.
The parameter coefficients in this model are the following:
\begin{align}
&\underline{\text{Coefficients}}&& \underline{\text{Name}}  
\nonumber \\ 
&\rho_b                                  && \text{Density\ of \ poroelastic \ matrix}; 
\nonumber
\\ 
&\lambda_b, \mu_b               && \text{Lam\'{e} \ constants};
\nonumber
\\ 
&\alpha_b                              &&\text{Biot-Willis \ coefficient \ \cite{show1}};
\label{parameters}
 \\ 
&c_b                                      && \text{Constrained \ storage\ coefficient};
\nonumber
\\ 
&k_b                                      && \text{Permeability \ of \ the\  poroelastic\  matrix.} 
\nonumber
\end{align}
In our work,  we will allow the permeability  $k_b$ of the poroelastic matrix to be a nonlinear function of one variable, which we specify below in Assumption \eqref{assumption_permeability}. We will also consider the case when the permeability $k_b$ is a given function of $\mathbf x$ and $t$. {All other coefficients listed above will be taken constant in this analysis.}

System \eqref{Biot_main} is obtained from the following system describing the balance of linear momentum
for the fluid-solid mixture, and conservation 
of mass for the fluid content, denoted by $\zeta$ (see \eqref{fluid_content} below to recall the definition of fluid content):
\begin{equation}\label{Biot_basic}
\begin{cases}
\rho_b\bfeta_{tt} -\nabla\cdot\mathbf \bfsigma_b = \mathbf F_b(\mathbf x, t),  &\hskip1cm   {\rm in} \ \Omega_b \times (0,T), \\[.2cm] 
\zeta_{t} + \nabla\cdot\mathbf{u}_b=S(\mathbf x, t), &\hskip1cm   {\rm in} \ \Omega_b \times (0,T),
\end{cases}
\end{equation}
where $\mathbf{u}_b$ denotes the discharge/filtration velocity associated with the pore pressure $p_b$, and $ \bfsigma_b$ denotes the
poroelastic stress tensor, which is given by:
$$
 \bfsigma_b = \bfsigma^E - \alpha_b p_b \boldsymbol{I},\hskip1cm   {\rm in} \ \Omega_b \times (0,T),
$$
where $\bfsigma^E$ denotes a given elasticity stress tensor, and $\boldsymbol{I}$ represents the identity tensor.
To close the system, two constitutive laws need to be prescribed, plus Darcy's law:
\begin{enumerate}
\item The constitutive law $\bfsigma^E = \bfsigma^E(\bfeta)$, describing the elastic material properties:
we assume the Saint Venant-Kichhoff material: 
\begin{align}
\boldsymbol\sigma^E =& 2\mu_b \boldsymbol{D}(\boldsymbol\eta) + \lambda_b {\rm tr} \boldsymbol{D}(\boldsymbol\eta) \boldsymbol{I}
=2\mu_b \boldsymbol{D}(\boldsymbol\eta) + \lambda_b \nabla \cdot \boldsymbol\eta \boldsymbol{I},
\end{align}
where  $\lambda_b,\mu_b$ denote the Lam\'{e} constants, and $\boldsymbol{D}(\boldsymbol\eta) = \frac{1}{2}(\nabla \boldsymbol\eta + (\nabla \boldsymbol\eta)^T)$ represents the symmetrized gradient. 
\item The constitutive law $\zeta = \zeta(p_b,\bfeta)$, defining the fluid content in terms of 
the fluid pore pressure, and the ``compressibility of the solid
matrix'', i.e., the change in the pore volume:
\begin{equation}\label{eq:fluid_increment}
 \zeta =  c_bp_b+\alpha_b\nabla \cdot \bfeta,
\end{equation}
 where $c_b\ge 0$ and $\alpha_b>0$ were introduced in \eqref{parameters}. Note that the case $c_b=0$ is permitted below, and this corresponds to {\em incompressible constituents} in the poroelastic material.
\item
 Darcy's law:
\begin{equation}
\label{eq:vel_darcy}
\mathbf{u}_b = - {\boldsymbol K}\nabla  p_b, 
\quad \mbox{with} \quad
{\boldsymbol K} = k_b \boldsymbol{I},
\end{equation}
where $\boldsymbol{K}$ is the permeability tensor.
\end{enumerate}

By using these constitutive laws in system \eqref{Biot_basic}, one recovers \eqref{Biot_main}. 
Notice that system \eqref{Biot_main} is obtained under the assumption of small deformations.

\begin{assumption}\label{assumption_permeability}
When considering the case of nonlinear permeability, we assume that the permeability tensor $\boldsymbol K$ is a {\sl{nonlinear function}} of {\sl{porosity}} $\phi$:
\begin{equation}
{\boldsymbol K} = k_b \boldsymbol{I}, \ {\rm with} \ k_b=k_{\text{ref}} f_k(\phi),
\end{equation}
where $k_{\text{ref}}$ is a reference value for the permeability of the mixture given in terms
of the volumetric fraction $\phi$ of the fluid component, called {\sl{porosity}}.
The particular form of the relationship between the permeability $k_b$ and the porosity $\phi$ is represented by the function $f_k(\phi)$, and it  
depends on the geometrical architecture of the pores inside the matrix and the physical properties of the fluid \cite{Hsu,causin2014poroelastic}.
\end{assumption}

For completeness, we recall the definitions of {\sl{porosity}} $\phi$ and {\sl{fluid content}} $\zeta$.
Porosity $\phi$ is defined as the volumetric fraction 
\begin{equation}\label{porosity}
\phi = \frac{V_f(\mathbf{x},t)}{V(\mathbf{x},t)}
\end{equation}
of the volume $V_f(\mathbf{x},t)$ occupied by the fluid within the poroelastic medium, and the representative elementary 
volume of the poroelastic medium $V(\mathbf{x},t)$, centered at $\mathbf{x} \in \Omega_b$ at time $t$. 

Similarly, if we denote by $V_s(\mathbf{x},t)$ the volume occupied by the solid, then under the assumption 
of fully saturated mixture, 
 the volumetric fraction 
${V_s(\mathbf{x},t)}/{V(\mathbf{x},t)} = 1 - \phi(\mathbf{x},t)$.

{\sl{Fluid content}} $\zeta$ is the increment in the volumetric fraction of the fluid component $\phi$ 
with respect to its baseline value $\phi_0$:
\begin{equation}\label{fluid_content}
\zeta(\mathbf{x},t) = \phi(\mathbf{x},t)-\phi_0(\mathbf{x}).
\end{equation}
Notice that the constitutive law for the fluid content \eqref{eq:fluid_increment} implies that porosity is given by
$$
 \phi = \phi(p_b,\bfeta) = \phi_0 + c_bp_b+\alpha_b\nabla \cdot \bfeta.
$$ 
 Thus, under the nonlinearity assumption, specified in Assumption~\ref{assumption_permeability}, we have:
 \begin{equation}\label{nonperm} 
 k_b=k_b(\phi)=k_b(\phi(p_b,\nabla \cdot \bfeta)) = k_b(c_bp_b+\alpha_b\nabla \cdot  \bfeta).
 \end{equation} 
For analysis purposes, we will be assuming the following properties of the permeability function $k_b$: 
\begin{assumption}
We assume that the permeability function $k_b: \reals \to \reals$ is continuous and that there exist constants $k_{\rm min} > 0$ and $k_{\rm max} >0$ s.t. 
$$0 < k_{\rm min} \leq k_b(\zeta) \leq k_{\rm max}, \ \ \forall \zeta \in \reals.$$
\end{assumption}
In considerations of the {\em linear} dynamics below, we will consider $k_b=k_b(\mathbf x, t)$ to be a given space and time dependent function. In the nonlinear considerations,
we will denote by
$k_b(\Psi(x))$ the Nemytskii operator associated with $k_b$. 
In that situation, our assumptions on the function $k_b$, and the theory of superposition operators \cite{RR, Tr}, 
will yield  that the operator $k_b$ is bounded and continuous from $L^2(\Omega \times (0,T))$ into $L^2(\Omega \times (0,T))$. In our considerations of uniqueness for the nonlinear problem below, we will consider the hypotheses that $k_b$ is a globally Lipschitz function on $\mathbb R$. 

\begin{remark} As in \cite{bgsw}, we could also consider a poroviscoelastic material, 
for which the stress tensor $\bfsigma_b$ would include a viscoelastic component, for example:
\begin{align} 
& \bfsigma_b = \bfsigma^E+\bfsigma^V - \alpha_b p_b \mathbf{I},
\qquad
\boldsymbol\sigma^V = 2\mu_V \boldsymbol{D}(\boldsymbol\eta_t) + \lambda_V {\rm tr} \boldsymbol{D}(\boldsymbol\eta_t) \boldsymbol{I},&
\end{align} 
where $\lambda_V$ and $\mu_V$ are independent Lam\'{e} parameters for viscoelasticity.
All the results from this manuscript would hold for the viscoporoelastic case as well.
\end{remark}


\if 1 = 0
TODO: Additional assumptions to be mentioned below:
{\em negligible inertia, 
small deformations and intrinsic incompressibility of each mixture component} ({\Lo These assumptions should be replaced just with "small deformations" here. Below we can mention that we include the case of negligible inertia and incompressible constituents})
\fi


\subsubsection{Poroelastic Plate Model in $\Omega_p$}

Domain $\Omega_p$ denotes a thin, poroelastic plate, separating the fluid flow region $\Omega_f$, 
from the thick, poroelastic structure in $\Omega_b$. 
A plate is a thin, 3D region,  bounded by two surfaces of small curvature,
whose distance defines the thickness of the plate. 
In our case, the plate will be assumed to be flat and of uniform thickness $h$, with the two surfaces located
at $h/2$ distance on either side of the {\sl{middle surface}} of the plate, denoted by $\omega_p$. 
The reference configuration of the middle surface $\omega_p$ is $x_3 = 0$.

We will be assuming that the elastic plate is porous, isotropic, and saturated by a viscous fluid. 
To capture the elastic deformation of the poroelastic skeleton, as well as fluid flow through the skeleton pores,  
we will be using a model which is a dynamic version of the  quasi-static Biot poroelastic plate model, 
studied by Biot in \cite{Biot55}, and rigorously justified 
by Marciniak-Czochra and Mikeli\'{c} in \cite{mikelic}. 
This model was obtained under the following hypotheses.
\begin{assumption}\label{KirchhoffPlate}
{{Poroelastic plate hypotheses:}}
\begin{enumerate}
\item The Kirchhoff hypothesis: Every straight line in the plate that was originally
perpendicular to the plate's middle surface, remains straight and perpendicular
to the deflected middle surface after deformation;
\item The fluid velocity derivatives in the longitudinal and lateral directions
are small compared to the transverse one.
\end{enumerate}
\end{assumption}
Therefore, the dominant flow is in the transverse direction.

The resulting equations age given in terms of displacement $\boldsymbol{w}=\boldsymbol{w}(x_1,x_2)$ 
of the plate's middle surface from its reference configuration $\omega_p$, where $(x_1,x_2)\in\omega_p$,
and the fluid pore pressure $p_p$, where $p_p$ is a function of three variables:
two corresponding to the coordinates along the middle surface $\omega_p$,
and the third one, denoted by $s \in (-h/2, h/2)$ corresponding to the {\sl{local transverse}} coordinate, so that
$p_p=p_p(x_1,x_2,s)$, 
where $(x_1,x_2,s) \in \Omega_p$.


According to \cite{mikelic}, for small deflections and $s$-driven pressure gradients, 
the in-plane and transverse plate dynamics fully decouple. 
Since in our problem the stress in the transverse direction is
larger than the stress in the longitudinal and lateral directions,
we consider only the {\em{transverse plate displacement}}, 
which we denote by $w(x_1,x_2,t)$, where  $(x_1,x_2) \in \omega_p$. 
For each $s$ slice in $(-h/2,h/2)$, the $w$ dynamics is constant. 
However, the fluid pore pressure $p_p(x_1,x_2,s)$ can vary for $s \in (-h/2,h/2)$.
The resulting reduced Biot model for poroelastic plate is then given by two equations,
one describing poroelastodynamics, which will be defined in $\omega_p$, 
and one describing conservation of mass of the fluid phase, given in terms of the evolution of $p_p$,
defined in $\Omega_p$ \cite{Biot55,mikelic,Fung,cowin1999bone}:
\begin{equation}
\begin{cases}
\rho_p w_{tt}+D\Delta^2_{\omega_p}w+\gamma w+\alpha_p\Delta_{\omega_p} \ds \int_{-h/2}^{h/2}sp_p~ds=
F_p(x_1,x_2,t)~~\text{ in }~~\omega_p \times (0,T),\\[.2cm]
\partial_t[c_pp_p-\alpha_p s \Delta_{\omega_p} w]- \partial_{s}(k_p \partial_{s} p_p)=0~~\text{ in }~~\Omega_p \times (0,T).
\end{cases}
\end{equation}
All Laplacians above refer to the in-plane Laplacian, i.e., $\Delta_{\omega_p}=\partial_{x_1}^2+\partial_{x_2}^2.$ 
The constant $D >0$ is the elastic stiffness coefficient for the plate, 
and $\gamma>0$ is an [elastic] coefficient that we added for technical reasons to provide coercivity,
since coercivity is not automatically satisfied in the periodic framework.
The constants $\rho_p, \alpha_p, c_p$ and $k_p$ are the same as those in 
\eqref{parameters}, taken for the poroelastic plate. One distinction for the poroelastic plate is that our results below will require $c_p>0$.
The source term $F_p$ corresponds to the loading of the poroelastic plate, which will come
from the jump in the normal components of the normal stress (traction) between the fluid on one side,
at the thick Biot poroelastic structure on the other: 
$F_p(x_1,x_2,t) ={{\bfsigma_b{\bf e}_3\cdot{\bf e}_3\big|_{\omega_p}-\bfsigma_f{\bf e}_3\cdot{\bf e}_3|_{\omega_p}} }$.
This will be specified below in the coupling conditions \eqref{dynamic}.
As before, Darcy's law, relating the filtration velocity to the pressure gradient, holds here:
\begin{equation}\label{DarcyPlate}
u_p= -k_p \partial_s p_p, \quad {\rm in} \quad \Omega_p \times (0,T),
\end{equation}
where $u_p$ denotes the transverse (dominant) component of filtration velocity, relative to the motion of the poroelastic matrix.

\subsubsection{Time-dependent Stokes Model in $\Omega_f$}
In the lower half space, denoted by $\Omega_f$, we consider the flow of a viscous, incompressible fluid
modeled by the time-dependent Stokes equations:
\begin{equation}\label{NS}
\left.
\begin{array}{l}
\displaystyle{\rho_f \frac{\partial \mathbf{u}}{\partial t} = \nabla \cdot \boldsymbol\sigma_f + \boldsymbol f} \\
\nabla \cdot \mathbf{u} = 0
\end{array}
\right\} \textrm{in}\; ~~\Omega_f \times (0,T),
\end{equation}
where $\mathbf{u}$ is the fluid velocity, $\boldsymbol\sigma_f=  2 \mu_f \boldsymbol{D}(\mathbf{u})-p_f \boldsymbol{I} $  is the 
Cauchy  stress tensor, $p_f$ is the fluid pressure, $\rho_f$ is the fluid density,  and $\mu_f$ is the fluid viscosity 
with  $\boldsymbol{D}(\mathbf{u})$ denoting the fluid strain-rate tensor (symmetrized gradient of velocity).

The fluid interacts with the multilayered poroelastic structure by exerting stress onto the poroelastic structure
that causes deformation of the poroelastic
skeleton, while at the same time generating the pressure difference, or more generally, the normal stress difference across 
the poroelastic plate, causing filtration flow through the poroelastic plate and on to the Biot poroelastic medium. 
The poroelastic plate serves as a fluid-structure interface with mass, which will provide a regularizing mechanism 
in the analysis of this fluid-structure interaction problem. 

In summary, the following are the three models holding in the three different subdomains of $\Omega$:
\begin{align}\label{Biot_sys}
&\begin{cases}
\rho_b\bfeta_{tt}-\mu_b\Delta \bfeta-(\mu_b+\lambda_b)\nabla (\nabla \cdot \bfeta)+\alpha_b \nabla p_b=\mathbf F_b, &\hskip1cm  (x_1,x_2,x_3) \in \Omega_b, \\[.2cm]
[c_bp_b+\alpha_b\nabla \cdot \bfeta]_t-\nabla \cdot (k_b\nabla p_b)=S, &\hskip1cm  (x_1,x_2,x_3) \in \Omega_b,\\[.2cm]
\end{cases}\\
&\begin{cases}\label{plate_sys}
\rho_{p} w_{tt}+D\Delta^2_{\omega_p}w+\gamma w+\alpha_p\Delta_{\omega_p} \ds \int_{-h/2}^{h/2}sp_p~ds=F_p, & (x_1,x_2) \in  \omega_p, 
\\[.2cm]
[c_pp_p-\alpha_ps\Delta_{\omega_p} w]_t- \partial_{s}(k_p \partial_{s} p_p)=0, &(x_1,x_2,s) \in \Omega_p,  \\[.2cm]
\end{cases}\\
&\begin{cases}\label{Stokes_sub}
\rho_f\mathbf u_t-\mu_f\Delta \bu+\nabla p_f = \mathbf f, &\hskip1cm  (x_1,x_2,x_3) \in \Omega_f,\\
\nabla \cdot \bu =0,  &\hskip1cm  (x_1,x_2,x_3) \in \Omega_f.
\end{cases}
\end{align}
\begin{remark} The notation for $k_b$ here (and below) in the problem description is general, so that it can accommodate the linear and nonlinear cases simultaneously. \end{remark}

\subsubsection{Coupling Conditions} 
The coupling between the fluid and multilayered poroelastic structure is assumed at the poroelastic plate,
which serves as a fluid-structure interface with mass.
The coupling conditions describe the relationship between the kinematic quantities, such as velocitites (and consequently
displacements), and the balance of forces. 
{{One set of coupling conditions will be prescribed across the middle surface of the plate $\omega_p$, and the other set of coupling conditions will be prescribed across the top and bottom surfaces on the plate, which we denoted by
$\omega_p^+$ and $\omega_p^-$. The top and bottom surfaces are defined in terms of the ``inflation variable'' $s$ as follows:
\begin{equation}\label{inflated_surfs}
\omega_p^+ = (0,1)^2 \times \{h/2\}=\{(x_1,x_2,s) | (x_1,x_2) \in \mathbb{R}^2, s=1/2\}
, \quad {\rm and} \quad \omega_p^- = (0,1)^2 \times \{-h/2\}.
\end{equation}

The coupling conditions along the middle surface of the plate $\omega_p$ are in line with the Kirchhoff hypothesis for plates,
see Assumption~\ref{KirchhoffPlate}, in the sense that those conditions are given in terms of 
 displacement $w$ of the middle surface of the plate, which is only a function of $x_1$ and $x_2$,
 therefore the coupling is along $\omega_p$. 
Across $\omega_p$, we prescribe the behavior of the {\emph{structure}} velocity $\partial_t w$ 
(the kinematic coupling condition), 
and the forces that drive the motion of the middle surface of the plate (the dynamic coupling condition).
Across $\omega_p^+$ and $\omega_p^-$ we prescribe continuity of vertical components 
of {\emph{fluid}} velocity (the kinematic coupling condition) and balance of normal component of normal stress (the dynamic coupling condition).

One novelty of the proposed model is the coupling in the pressure. See Remark~\ref{remark:pressure}.
As we shall see below, the coupling of the plate filtration pressure at the top interface $\omega_p^+$ and at the bottom interface $\omega_p^-$
with the surrounding fluid pressure, will give rise to a pressure jump between the top and bottom interface,
which is what drives the filtration flow through the plate.
}}

{\em 1. Coupling across the plate's middle surface $\omega_p$:} We have two sets of coupling conditions,
the kinematic and dynamic.
For the kinematic coupling condition, we assume continuity 
between the plate displacement and the Biot displacement evaluated at the
interface $\omega_p$, 
while on the fluid side we assume the Beavers-Joseph-Saffman coupling condition, stating
that the tangential components of the fluid velocity evaluated at $\omega_p$ are proportional 
to the corresponding tangential components of the fluid normal stress evaluated at $\omega_p$.
The constant of proportionality, denoted by $\beta$, denotes the {\em{slip length}}.
For the dynamic coupling, see \eqref{dynamic} below, we state that the elastodynamics of the poroelastic plate
is driven by 
the jump in the normal components of normal stress between the fluid and Biot poroelastic medium,
where the jump is evaluated at the two-dimensional fluid-structure interface $\omega_p$.
Here $F_p$ is the forcing from the plate equation \eqref{plate_sys}:
\begin{align}
&\left\{
\begin{array}{rcl}
\langle 0,0,w\rangle &=& \bfeta \big|_{\omega_p},\\
\beta\mathbf u \big|_{\omega_p}\cdot \mathbf e_i &=& {-} [\boldsymbol{\sigma}_f\big|_{\omega_p} \mathbf e_3]\cdot\mathbf e_i,~~i=1,2,
\end{array}
\right.
\label{kinematic}
\\
& \qquad \qquad \quad \ \ F_p = {{\bfsigma_b{\bf e}_3\cdot{\bf e}_3\big|_{\omega_p}-\bfsigma_f{\bf e}_3\cdot{\bf e}_3|_{\omega_p}} }.
 \label{dynamic}
\end{align}
{{
We remark here that the first condition in \eqref{kinematic} in the continuity of vertical components 
of structure velocities (assuming equal structure displacements initially), and the second condition in \eqref{kinematic} in fact reads:
$$
(\beta\mathbf u \big|_{\omega_p} -  \partial_t \boldsymbol{w} )\cdot \mathbf e_i = {-} [\boldsymbol{\sigma}_f\big|_{\omega_p} \mathbf e_3]\cdot\mathbf e_i,~~i=1,2,
$$
where the left hand-side is the tangential velocity slip. However, since the full Lagrangian plate displacement $\boldsymbol{w} = \langle 0,0,w\rangle $, the  tangential components ($i=1,2$) of $\partial_t \boldsymbol{w}$ are both equal to zero.
}}

{\em 2. Coupling across the Biot-plate interface ($\omega_p^+$):} Here we specify the coupling
that involves the quantities that are defined on the ``inflated'' plate domain $\Omega_p$, such as filtration velocity 
and pore pressure, and couple them to the corresponding quantities in the Biot model. 
For this purpose, we evaluate the plate quantities at the inflated surface $\omega_p^+$, defined in 
\eqref{inflated_surfs}, and the Biot quantities at the interface portion of the boundary of the Biot domain, namely on $\omega_p$.
In particular, we assume continuity of filtration velocities, which are defined
relative to the motion of the poroelastic matrix: $\mathbf{u}_b \cdot {\bf e}_3 = u_p$, and continuity of pore pressures.
Using Darcy's laws \eqref{eq:vel_darcy} and \eqref{DarcyPlate}, the two coupling conditions read:
\begin{align}
&k_p\partial_{s}p_p\big|_{\omega_p^+}= k_b \partial_{x_3}p_b\big|_{\omega_p}& 
\label{velocity_omega+}\\
&p_p\big|_{\omega_p^+}=p_b\big|_{\omega_p}.&
\label{pressure_omega+}
\end{align}
{{
We remark here that condition \eqref{velocity_omega+} describes continuity of vertical components of filtration velocities
(relative to the elastic structure motion),
after accounting for Darcy's law. 
}}

{\em 3. Coupling across the Stokes-plate interface ($\omega_p^-$):} Similarly to the coupling presented above, 
here we specify the coupling of the quantities defined on the ``inflated'' plate domain $\Omega_p$, namely 
the plate filtration velocity and pressure, with the Stokes flow.  
Again, we evaluate the plate quantities at the inflated surface $\omega_p^-$ adjacent to the fluid domain, 
and the Stokes quantities at the interface boundary of the fluid domain, namely on $\omega_p$.
In particular, we assume continuity of filtration velocities relative to the motion of
the poroelastic matrix: $w_t-(\bu\cdot \mathbf e_3)\big|_{\omega_p} = u_p$, and continuity of 
the normal component of the normal fluid stress.
Using Darcy's law \eqref{DarcyPlate}, the two coupling conditions read:
\begin{align}\label{finalcouple}
&w_t-(\bu\cdot \mathbf e_3)\big|_{\omega_p}=k_p\partial_{s}p_p\big|_{\omega_p^-},&\\[.2cm]
&-\boldsymbol{\sigma}_{f}\mathbf e_3 \cdot \mathbf e_3\big|_{\omega_p}= p_p\big|_{\omega_p^-}.&
\label{pressure_omega-}
\end{align}
{{
We remark that condition \eqref{finalcouple} describes continuity of vertical components of  velocities (relative
to the structure motion),
where the fluid velocity on the left hand-side is corrected by the velocity of the structure $w_t$.
}}

\begin{remark}\label{remark:pressure}
By combining the coupling conditions for the fluid pressure given in \eqref{pressure_omega+} and \eqref{pressure_omega-},
we see that across the FPSI interface, the pressure (the normal component of normal stress, in fact)  jumps from 
${-}\boldsymbol{\sigma}_{f}\mathbf e_3 \cdot \mathbf e_3\big|_{\omega_p}$ on the fluid side
 to  $p_b\big|_{\omega_p}$ on the Biot side. 
In between,  a pressure gradient-driven vertical/transverse flow through the poroelastic plate takes place
from $p_p\big|_{\omega_p^-}$ to $p_p\big|_{\omega_p^+}$. 
Thus, the jump in the pressure at the interface $\omega_p$ 
(namely, the jump in the normal components of normal stress), given by
$
[p_b {-}\boldsymbol{\sigma}_{f}\mathbf e_3 \cdot \mathbf e_3]
$
on $\omega_p$, drives the filtration flow through the interface, i.e., through the poroelastic plate.

\if 1 = 0
After a calculation, one gets that the pressure jump is proportional to 
$\partial {\rm{u}}_{x_3}/\partial x_3$,
namely, that across the interface $\omega_p$, the following jump in the pressure holds:
$$
p_b - p_f = -\mu_f  \frac{\partial \rm{u}_{x_3}}{\partial x_3}.
$$
It is interesting to notice that a calculation based on boundary layer and homogenization 
approaches of the effective pressure across the interface between an incompressible, viscous fluid and a poroelastic medium,
obtained in \cite{JagMik96,JagMik00}, indicates that indeed, there should be a pressure jump across the interface.
However, Marciniak-Czochra and Mikeli\'{c} found in \cite{JagMik96,JagMik00} that the pressure jump should be proportional to the fluid shear stress,
namely that across the fluid-Biot poroelastic medium interface:
$$
p_b - p_f = C_{\text{pressure}} \frac{\partial \rm{u}_{\tau}}{\partial x_3}, 
$$
where $ \rm{u}_\tau$ is the tangential component of the fluid velocity, and $C_{\text{pressure}}$ is
determined from the boundary layer pressure, and is of order ${\cal{O}}(1)$ for non-isotropic porous media, for example. 
\fi
\end{remark}

\subsubsection{Boundary and Initial Conditions}\label{BCIC}

{\bf Boundary Conditions on $\partial\Omega$:} We take periodic boundary conditions for all quantities of interest in the longitudinal and lateral directions,
and identify $x_{i}=0$ with $x_i=1$, for $i = 1,2$.

 For $x_3=-1$, we consider the standard {\em{no-slip}} boundary conditions for the Stokes fluid, 
 i.e., $\mathbf u \equiv \mathbf 0\big|_{x_3=-1}$. 
For $x_3=1$, we take Neumann-type (homogeneous) boundary conditions for the Biot system. 

We summarize  the boundary conditions in the following table:

\begin{align}
\underline{\text{Region}}&& \underline{\text{Functions}} && \underline{\text{Boundary Condition}} \\ 
\label{PeriodicBC}
x_i=0, ~x_i=1, i= 1,2 && \boldsymbol\eta, w, p_b, \mathbf u && \text{periodic} \\
x_3=-1 && \mathbf u && \mathbf u\big|_{x_3=-1}= \mathbf 0\\
x_3=1 && \boldsymbol \eta,~p_b &&\boldsymbol \sigma(\boldsymbol \eta,p_b) \mathbf e_3 =\mathbf 0, ~\partial_{x_3}p_b=0
\label{zero_traction}
\end{align}
\begin{remark}{Note that we do not prescribe boundary conditions for $p_p$ on the lateral boundaries, owing the the structure of the pressure equation for the poro-plate. This is to say, there is no differential operator acting on $p_p$ in the tangential directions, and thus no boundary conditions to be prescribed. We will see this also in the functions spaces prescribed below in Section \ref{secweak}.}\end{remark}

\noindent
{\bf Initial Conditions at $t = 0$:}
\begin{align}
{\rm Fluid}:&~\; \bu(.,0)=\bu^0, \label{FluidIC}
\\
{\rm Poroelastic\; plate:}&~\; w(.,0)=w^0,\; \partial_t w(.,0)=\dot w^0,\; p_p(.,0)=p_p^0,\label{PlateIC}
\\
{\rm Biot\; equations:}& ~\; \bfeta(.,0)=\bfeta^0,\; \partial_t\bfeta(,.0)=\dot{\bfeta}^0,\; p_b(.,0)=p_b^0.
\label{BiotIC}
\end{align}
Here, we used $\dot w^0$ and $\dot{\bfeta}^0$ to denote the initial data for the time derivative of $w$ and $\bfeta$, respectively.

\subsection{Quasi-static Problem}\label{sec:quasistatic*}

The {\em quasi-static} version of problem \eqref{Biot_sys}-\eqref{BiotIC} refers to the problem
in which the inertial effects in the Biot model \eqref{Biot_main} are neglected:

\begin{assumption}\label{quasistatic} The quasi-static model refers to setting $\rho_b = 0$ in \eqref{Biot_sys},
and allowing $\rho_p\geq 0$ in \eqref{plate_sys}.
\end{assumption}

This quasi-static assumption arises naturally in the classical Biot model of consolidation for
a linearly elastic and porous solid which is saturated by a slightly compressible
viscous fluid. See, e.g., \cite{show1}.
The quasi-static assumption requires modifying the initial conditions, but leaves the coupling conditions the same. 
As we shall see in Sec.~\ref{sec:quasistatic}, we can extend our linear existence result to the nonlinear case {\em when} Assumption \ref{quasistatic} is taken. On the other hand, it has been noted by many authors that the elliptic-parabolic nature of the quasi-static Biot problem presents its own challenges \cite{BW,show1}.

The new {\em{initial conditions}}, owing to loss of acceleration terms, are given by:
\begin{itemize}
\item If $\rho_b=0$, initial conditions \eqref{BiotIC} are replaced with initial condition for the fluid content:
\begin{align}\label{QBiotIC}
\zeta_b(.,0)=(c_bp_b+\alpha_b\nabla\cdot\bfeta)(.,0)=\zeta_b^0.
\end{align}
\item If $\rho_p=0$, initial conditions \eqref{PlateIC} are replaced with initial condition for the poroelastic plate fluid content:
\begin{align}\label{QPlateIC}
\zeta_p(.,0)=(c_pp_p-\alpha_p s \Delta_{\omega_p} w)(.,0)=\zeta_p^0.
\end{align}
\end{itemize}
In the analysis of such quasi-static Biot models (see \cite{BW, zenisek, bgsw, cao, show1,owc}), the structure of the initial condition is a non-trivial issue.  In this treatment, we only consider $\zeta^0_b$ for which there exist $\bfeta^0 = \bfeta(0)$ and $p^0_b = p_b(0)$ s.t. $c_bp^0_b+\alpha_b\nabla\cdot\bfeta^0 = \zeta^0_b$. Such a hypothesis is common \cite{cao,zenisek,bgsw}.
\begin{remark}
For a very weak notion of solution in the linear case, the work in \cite{showmono,show1} only requires initial conditions as specified above, which is to say, for the respective fluid contents alone. However, in nearly every other construction, additional regularity of the initial condition is required. This is a by-product of using a priori estimates to construct solutions. See \cite[Remark 2]{zenisek} and \cite[p.327]{show1} where this is explicitly discussed.
\end{remark}

Details of the existence result are given in Section~\ref{sec:quasistatic}.



\section{Energy, Function Spaces, and Weak Formulation}
In this section we derive a formal energy estimate, present functions spaces for the solution, and define weak solution.
For this purpose, we recall 
the principal solutions variables in each sub-region:
\begin{itemize}
\item $(\bu,p_f)$---Stokes fluid velocity and pressure defined on $\Omega_f$;
\item $(w,w_t,p_p)$---transverse displacement and velocity of the poroelastic plate and poroelastic plate pressure, defined on $\omega_p$ and $\Omega_p$ respectively; 
\item $(\bfeta,\bfeta_t,p_b)$---poroelastic displacement, velocity and pressure defined on $\Omega_b$ (thick Biot region). The quantity $\bu_b=-\boldsymbol K \nabla p_b=-k_b(c_bp_b+\alpha_b\nabla \cdot \bfeta)\nabla p_b$ is the discharge velocity in thick poroelastic material.
\end{itemize} 

\subsection{A Formal Energy Inequality}
By testing equations \eqref{Biot_sys}-\eqref{Stokes_sub} by the solution variables, specified above,
and by using coupling conditions \eqref{kinematic}-\eqref{pressure_omega-}, and periodic boundary conditions, we obtain the following formal energy identity: 


\begin{align}\label{Energy}
&\phantom{+} \ \frac{1}{2}\frac{d}{dt}\Big (
\rho_f\|\bu\|^2_{L^2(\Omega_f)}
+\rho_p \|w_t\|^2_{L^2(\omega_p)}+\rho_b\|\partial_t\bfeta\|^2_{L^2(\Omega_b)} \Big)&  
\\
&+\dfrac{1}{2}\dfrac{d}{dt}\Big(c_p\|p_p\|^2_{L^2(\Omega_p)}
+\|\bfeta\|^2_{E}+c_b\|p_b\|^2_{L^2(\Omega_b)}+D\|\Delta w\|^2_{L^2(\omega_p)}+{ \gamma}\| w\|^2_{L^2(\omega_p)}
\Big )&
\nonumber
\\
&+{2}\mu_f\|\boldsymbol{D}(\bu)\|^2_{L^2(\Omega_f)}+ \beta\|\bu\cdot\boldsymbol\tau\|^2_{L^2({\omega_p})}+
\|k_p^{1/2}\partial_{s} p_p\|^2_{L^2(\Omega_p)}+\|k_b^{1/2}\nabla p_b\|^2_{L^2(\Omega_b)}\leq C_{{data}},&
\nonumber
\end{align}
where 
\begin{equation*}\label{Enorm}
\|\bfeta\|^2_E=\ds \int_{\Omega_b}[\bfsigma^E:\nabla\bfeta] d\mathbf x =2\mu_E(\boldsymbol D(\bfeta),\boldsymbol D(\bfeta))_{L^2(\Omega_b)}+\lambda_E (\nabla \cdot \bfeta, \nabla \cdot \bfeta)_{L^2(\Omega_b)},
\end{equation*}
and 
$$||\bu \cdot \bt||_{L^2(\omega_p)}^2 = \sum_{i=1,2} ||\bu \cdot \bt_i||^2,$$ 
where $\bt_i$ are the tangent vectors,
given by $\mathbf e_1,~\mathbf e_2$ on $\omega_p$ in this case.
Note that the formal energy identity for the quasi-static dynamics can be obtained by letting $\rho_b=0$, and this is valid for $\rho_p=c_b=0$ as well. Our results will include all these cases. {See Section \ref{genweak} and Remark \ref{testremark}.}
\begin{remark}
We shall see later that this energy inequality will hold for constructed solutions, but need not hold in the general case for weak solutions, even in the case of linear dynamics.
\end{remark}

\subsection{Function Spaces for Weak Solutions}\label{secweak}
For any function space $X$,  {let us introduce the notation that} $X_{\#}$ is the space of all functions from $X$ that are $1$-periodic in directions $x_1$ and $x_2$:
\begin{align*}
X_{\#}:=\{f\in X: f\big|_{x_1=0}=f\big|_{x_1=1},\;f\big|_{x_2=0}=f\big|_{x_2=1}\}.
\end{align*}
Stokes equation:
\begin{align*}
H=&\{\bu\in \mathbf L^2(\Omega_f): \nabla\cdot\bu=0,\; (\bu\cdot{\bf e}_3)\big |_{x_3=-1}=0\},
\\
V=&\{\bu\in \mathbf H^1(\Omega_f): \nabla\cdot\bu=0,\; \bu\big |_{x_3=-1}=0\},
\\
\mathcal{V}_f=&L^{\infty}(0,T;H)\cap L^2(0,T;V_{\#}).
\end{align*}
Poroelastic plate:
\begin{align*}
H^{0,0,1}=&\{p\in L^2(\Omega_p):\partial_{s} p\in L^2(\Omega_p)\}
\\
\mathcal{V}_p=&\{w\in L^{\infty}(0,T;H^2_{\#}(\omega_p)): \rho_p w\in W^{1,\infty}(0,T;L^2(\omega_p))\},
\\[.2cm]
\mathcal{Q}_p=&\big\{p\in L^2(0,T;H^{0,0,1}) : c_p p\in L^{\infty}(0,T;L^2(\Omega_p))\big\}.
\end{align*}
Biot equations:
\begin{align*}
\mathcal{V}_b=&\{\bfeta\in L^{\infty}(0,T;\mathbf H^1_{\#}(\Omega_b)):\rho_b\bfeta\in W^{1,\infty}(0,T;\mathbf L^2(\Omega_b))\},
\\[.2cm]
\mathcal{Q}_b=& \{p\in L^2(0,T;H^1_{\#}(\Omega_b)):c_b p\in L^{\infty}(0,T;L^2(\Omega_b))\}.
\end{align*}
Weak {\em solution space}:
\begin{equation*}\label{Solspace}
\mathcal{V}_{\rm sol}=\big\{(\bu,w,p_p,\bfeta,p_b)\in \mathcal{V}_f\times\mathcal{V}_p\times\mathcal{Q}_p\times\mathcal{V}_b\times\mathcal{Q}_b~:~{p_p\big |_{s=h/2}=p_b\big |_{x_3=0},\; \bfeta\big |_{x_3=0}=w{\bf e}_3}\big\}.
\end{equation*}

\begin{remark} We included the ``inertial" and ``compressibility" constants in the definition of the solution spaces in order to unify the exposition and, at the same time, point to the different regularity of solution for the different scenarios. \end{remark}

In conjunction with $\mathcal V_{\rm sol}$ we define the {\em test space} as follows:
\begin{align}\label{Testspace}
\mathcal{V}_{\rm test}=\big\{(\bv,z,q_p,\bfpsi,q_b) \in &~ C^1_c\big ([0,T); V_{\#}\times H^2_{\#}(\omega_p) \times H^{0,0,1}\times \mathbf H^1_{\#}(\Omega_b)\times  H^1_{\#}(\Omega_b)\big)
\\ \nonumber
&: 
{q_p\big |_{s=h/2}=q_b\big |_{x_3=0},\; \bfpsi\big |_{x_3=0}=z{\bf e}_3}\big\}
\end{align}
where we use the standard notation $C^1_c$ for continuously differentiable functions with compact support. 

\subsection{General Weak Formulation}\label{genweak}
We aim to construct a function $(\bu, w, p_p, \bfeta, p_b)\in \mathcal{V}_{\rm sol}$, call a {\sl weak solution}, such that
for every smooth test function $(\bv,z,q_p,\bfpsi,q_b)\in\mathcal{V}_{\rm test}$, the {\sl{weak formulation}} 
specified below in Definition~\ref{DefWeakSol} holds. To formulate the weak formulation we introduce the following notation:
\begin{itemize}
\item $((\cdot, \cdot))_{\mathscr O}$ denotes an inner-product on $L^2(0,T,L^2(\mathscr O))$, for a domain $\mathscr O$;
\item $\langle \cdot, \cdot\rangle$ denotes a particular ``duality" pairing whose precise interpretation will be specified below.
\end{itemize}

To obtain the weak formulation for the coupled system, we proceed by formally multiplying the equations in \eqref{Biot_sys}-\eqref{Stokes_sub}
 by test functions and integrating by parts as described below. 

We begin by using the test function $\bfpsi  \in C^1_c \big([0,T); \mathbf H^1_{\#}(\Omega_b) \big)$, with $ \bfpsi\big |_{x_3=0}=z{\bf e}_3$ 
for the balance of linear momentum equation in the Biot model  \eqref{Biot_sys}, to obtain
\begin{equation}\label{Biot1_varform}
-\rho_b((\bfeta_t,\boldsymbol \psi_t))_{\Omega_b}+((\boldsymbol \sigma_b(\bfeta,p_b),\nabla \boldsymbol \psi))_{\Omega_b} -
 \int_0^T \int_{\partial \Omega_b} \boldsymbol \sigma_b(\bfeta,p_b) \mathbf n \cdot  \bfpsi 
 = \rho_b(\bfeta_t,\boldsymbol \psi)\big|_{t=0} +  ((\mathbf F_b,\boldsymbol \psi))_{\Omega_b}.
\end{equation} 
For the balance of mass equation, we use the test function $q_b \in C^1_c \big([0,T);H^1_{\#}(\Omega_b) \big)$, and obtain
\begin{align}\label{Biot2_varform}
-((c_bp_b+\alpha_b\nabla \cdot \bfeta , \partial_tq_b))_{\Omega_b} &+((k_b\nabla p_b,\nabla q_b))_{\Omega_b} - \int_0^T\int_{\partial \Omega_b} k_b\nabla p_b \cdot \mathbf n \ q_b
 \\ 
=& (c_bp_b+\alpha_b\nabla\cdot \boldsymbol \eta ,q_b)\big|_{t=0} +((S,q_b))_{\Omega_b}. 
\nonumber
\end{align}
For the equation describing the poroelastic plate dynamics \eqref{plate_sys}, we use test function $z \in C^1_c\big ([0,T); H^2_{\#}(\omega_p)\big)$ to obtain
\begin{align}\label{plate1_varform}
-\rho_p((w_t,\partial_tz))_{\omega_p}+D((\Delta w,& \Delta z))_{\omega_p}+\alpha_p((\int_{-h/2}^{h/2} ( sp_p) ds,\Delta z))_{\omega_p}+\gamma ((w,z))_{\omega_p}
\\ 
=
&~\rho_p(w_t,z)\big|_{t=0} + \int_0^T \int_{\omega_p}F_p z. 
\nonumber 
\end{align}
In the conservation of mass of the fluid phase for the poroelastic plate we use the test function $q_p \in C^1_c\big ([0,T); H^{0,0,1}_{\#}\big)$ to obtain
\begin{align}\label{plate2_varform}
-((c_pp_p-\alpha_ps\Delta w,\partial_t q_p))_{\Omega_p}
&+((k_p\partial_{s}p_p,\partial_{s}q_p))_{\Omega_p} 
- \int_0^T \int_{\partial \Omega_p} k_p\partial_sp_p n_3 q_p
\\ 
  =&~ (c_pp_p-\alpha_ps\Delta w,q_p)_{\Omega_p}\big|_{t=0}.
  \nonumber
\end{align}
Finally, we use test function $\bv \in C^1_c\big ([0,T);V_{\#}\big)$ in the Stokes equation \eqref{Stokes_sub}, which provides
\begin{equation}\label{Stokes_varform}
 -\rho_f((\bu,\bv_t))_{\Omega_f}+\mu_f ((\boldsymbol D(\bu),\boldsymbol D(\bv)))_{\Omega_f} + \int_0^T \int_{\partial \Omega_f} \boldsymbol \sigma_f \mathbf n \cdot \bv 
 = \rho_f(\bu,\bv)\big|_{t=0} + ((\mathbf f,\bv))_{\Omega_f}.
\end{equation}
Note that when we combine \eqref{Biot1_varform} - \eqref{Stokes_varform}, the collection of boundary terms is represented by
\begin{align} \label{boundaryterms_varform} 
I_{\rm bdry} =&~ - \int_0^T \int_{\partial \Omega_b} \boldsymbol \sigma_b(\bfeta,p_b) \mathbf n \cdot  \bfpsi  - \int_0^T \int_{\omega_p}\boldsymbol \sigma_b \mathbf e_3 \cdot \mathbf e_3 \big|_{\omega_p} z 
\\ 
&-\int_0^T\int_{\partial \Omega_b} k_b\nabla p_b \cdot \mathbf n \ q_b  - \int_0^T \int_{\partial \Omega_p} k_p\partial_sp_p n_3 q_p 
\nonumber
\\
&+  \int_0^T \int_{\omega_p}\boldsymbol \sigma_f \mathbf e_3 \cdot \mathbf e_3 \big|_{\omega_p} z  + \int_0^T \int_{\partial \Omega_f} \boldsymbol \sigma_f \mathbf n \cdot \bv. 
\nonumber
\end{align}
After using the boundary condition $\boldsymbol \sigma(\boldsymbol \eta,p_b) \mathbf e_3 =\mathbf 0$ and the fact that $ \bfpsi\big |_{x_3=0}=z{\bf e}_3$, 
the first term in \eqref{boundaryterms_varform} becomes
\begin{equation}\label{bdytermelastic}
- \int_0^T \int_{\partial \Omega_b} \boldsymbol \sigma_b(\bfeta,p_b) \mathbf n \cdot  \bfpsi = -\int_{x_3 =1} \boldsymbol \sigma_b\mathbf e_3 \cdot  \bfpsi  + \int_{\omega_p} \boldsymbol \sigma_b\mathbf e_3 \cdot  \bfpsi = \int_0^T \int_{\omega_p}\boldsymbol \sigma_b \mathbf e_3 \cdot \mathbf e_3 \big|_{\omega_p} z, 
\end{equation}
and therefore the first two terms cancel each other. 
The next two terms in \eqref{boundaryterms_varform} can be rewritten as follows
\begin{align}\label{terms3and4}
- \int_0^T\int_{\partial \Omega_b} &k_b\nabla p_b \cdot \mathbf n \ q_b   - \int_0^T \int_{\partial \Omega_p} k_p\partial_sp_p n_3 q_p
\\
=&~ - \int_{x_3=1} k_b \partial_{x_3}p_b q_b + \int_{\omega_p} k_b\partial_{x_3}p_bq_b + \int_{s =- h/2} k_p \partial_sp_p q_p - \int_{s = h/2}k_p \partial_sp_p q_p.
\nonumber
\end{align}
By using the boundary condition $\partial_{x_3}p_b = 0$ at $x_3 = 1$, 
the coupling across the Biot-plate interface $k_p\partial_{s}p_p\big|_{\omega_p^+}= k_b \partial_{s}p_b\big|_{\omega_p}$,
and the coupling across the Stokes-plate interface 
$w_t-(\bu\cdot \mathbf e_3)\big|_{\omega_p}=k_p\partial_{s}p_p\big|_{\omega_p^-}$, equation \eqref{terms3and4} simplifies to 
\begin{equation}\label{terms3and4-final}
- \int_0^T\int_{\partial \Omega_b} k_b\nabla p_b \cdot \mathbf n \ q_b  - \int_0^T \int_{\partial \Omega_p} k_p\partial_sp_p n_3 q_p
 =  \int_{s =- h/2} (w_t-(\bu\cdot \mathbf e_3)\big|_{\omega_p})q_p.
\end{equation}
Using the second coupling at the Stokes-plate interface $-\boldsymbol{\sigma}_{f}\mathbf e_3 \cdot \mathbf e_3\big|_{\omega_p}= p_p\big|_{\omega_p^-}$, we get 
\begin{equation}\label{bdyfluid_1}
\int_0^T \int_{\omega_p}\boldsymbol \sigma_f \mathbf e_3 \cdot \mathbf e_3 \big|_{\omega_p} z = - \int_0^T \int_{\omega_p}p_p\big|_{\omega_p^-} z.
\end{equation}
Lastly, using again the coupling at the Stokes-plate interface $-\boldsymbol{\sigma}_{f}\mathbf e_3 \cdot \mathbf e_3\big|_{\omega_p}= p_p\big|_{\omega_p^-}$, the boundary condition $\bv = 0$ on $\{x_3 = -1\}$, and the coupling across the plate's middle surface $\omega_p$,
$\beta\mathbf u \big|_{\omega_p}\cdot \mathbf e_i = {-} [\boldsymbol{\sigma}_f\big|_{\omega_p} \mathbf e_3]\cdot\mathbf e_i,~~i=1,2$, we can simplify the boundary term associated with the fluid as follows:
\begin{align} \label{bdyfluid_2} \int_0^T \int_{\partial \Omega_f} \boldsymbol \sigma_f \mathbf n \cdot \bv  =&~ \int_{x_3 = -1} \boldsymbol \sigma_f \mathbf e_3 \cdot \bv - \int_{\omega_p} \boldsymbol \sigma_f \mathbf e_3 \cdot \bv \\ \nonumber
 =&~ - \int_{\omega_p} [ \boldsymbol \sigma_f \mathbf e_3 \cdot \boldsymbol e_1 (\bv \cdot \boldsymbol e_1) + \boldsymbol \sigma_f \mathbf e_3 \cdot \boldsymbol e_2 (\bv \cdot \boldsymbol e_2) + \boldsymbol \sigma_f \mathbf e_3 \cdot \mathbf e_3 (\bv \cdot \mathbf e_3)]\\ \nonumber
=&~ \int_{\omega_p} \beta (\bu \cdot \boldsymbol t) (\bv \cdot \boldsymbol t) + \int_{\omega_p}p_p\big|_{\omega_p^-}( \bv \cdot \mathbf e_3).
\end{align}
Combining \eqref{boundaryterms_varform} with \eqref{bdytermelastic}, \eqref{terms3and4-final}, \eqref{bdyfluid_1}, and \eqref{bdyfluid_2}, we obtain that the sum of the boundary terms in the variational form simplifies to 
$$
I_{\rm bdry} = -((p_p\big|_{s=-h/2},z-\bv \cdot \mathbf e_3))_{\omega_p}+(\langle w_t-\bu \cdot \mathbf e_3,q_p\big|_{s=-h/2}\rangle)_{\omega_p} +\beta(( \bu \cdot \boldsymbol t , \bv \cdot \boldsymbol t ) )_{\omega_p}.
$$

We can now introduce the definition of weak solution for the fully coupled system.
The following weak form applies to all cases of admissible values of $\rho_b,\rho_p,c_b,c_p$.
Namely, the solution and test spaces are built to accommodate subsets of these parameters vanishing.

\begin{definition}\label{DefWeakSol}
We say that $(\bu, w, p_p, \bfeta, p_b)\in \mathcal{V}_{\rm sol}$ is a {\em weak solution} to \eqref{Biot_sys}--\eqref{BiotIC},
if for every test function 
$(\bv,z,q_p,\bfpsi,q_b)\in\mathcal{V}_{\rm test}$ the following identity holds:

\begin{align}
 \nonumber
&-\rho_b((\bfeta_t,\boldsymbol \psi_t))_{\Omega_b}+((\boldsymbol \sigma_b(\bfeta,p_b),\nabla \boldsymbol \psi))_{\Omega_b}-((c_bp_b+\alpha_b\nabla \cdot \bfeta , \partial_tq_b))_{\Omega_b}+((k_b\nabla p_b,\nabla q_b))_{\Omega_b}&\\[.2cm] 
\nonumber
&-((c_pp_p-\alpha_ps\Delta w,\partial_t q_p))_{\Omega_p}+((k_p\partial_{s}p_p,\partial_{s}q_p))_{\Omega_p} & 
\\[.1cm]
\nonumber
&-\rho_p((w_t,\partial_tz))_{\omega_p}+D((\Delta w,\Delta z))_{\omega_p}+\alpha_p((\int_{-h/2}^{h/2} [sp_p]ds,\Delta z))_{\omega_p}+\gamma ((w,z))_{\omega_p}&
\\
\label{WeakForm}
&-((p_p\big|_{s=-h/2},z-\bv \cdot \mathbf e_3))_{\omega_p}+(\langle w_t-\bu \cdot \mathbf e_3,q_p\big|_{s=-h/2}\rangle)_{\omega_p}& 
\end{align}
\begin{align}
\nonumber
& -\rho_f((\bu,\bv_t))_{\Omega_f}+\mu_f ((\boldsymbol D(\bu),\boldsymbol D(\bv)))_{\Omega_f}+\beta(( \bu \cdot \boldsymbol t , \bv \cdot \boldsymbol t ) )_{\omega_p}& 
\\[.3cm] 
\nonumber
=~&\rho_b(\bfeta_t,\boldsymbol \psi)\big|_{t=0}+(c_bp_b+\alpha_b\nabla\cdot \boldsymbol \eta ,q_b)\big|_{t=0}
+(c_pp_p-\alpha_ps\Delta w,q_p)_{\Omega_p}\big|_{t=0}+\rho_p(w_t,z)\big|_{t=0}+\rho_f(\bu,\bv)\big|_{t=0}&
\\[.2cm]
\nonumber
&+ ((\mathbf F_b,\boldsymbol \psi))_{\Omega_b}+((S,q_b))_{\Omega_b}+((\mathbf f,\bv))_{\Omega_f},&
\nonumber
\end{align}
where
\begin{align}\nonumber
\rho_b&\bfeta(\cdot,0)=\rho_b\bfeta^0,~\rho_b\bfeta_t(\cdot,0)=\rho_b\bfeta_t^0;\; \rho_p w(\cdot,0)=\rho_p w^0, ~\rho_pw_t(\cdot,0)=\rho_pw_t^0, 
\\[.1cm]
\nonumber
c_b&p_b(\cdot,0)+\alpha_b\nabla \cdot \bfeta(\cdot,0)=\zeta_b^0; \; c_pp_p(\cdot,0)-\alpha_p\Delta w(\cdot,0)=\zeta_p^0.&
\end{align}
\end{definition}
The interpretation of the 
dual mapping in term $\langle w_t-\bu \cdot \mathbf e_3,q_p\big|_{s=-h/2}\rangle_{\omega_p}$
is the following:
\begin{itemize}
\item If $\rho_p>0$, the dual mapping is just $L^2(0,T;L^2(\omega_p))$ scalar product, which is well defined since by $w_t\in L^2(0,T;L^2(\omega_p))$ (see definition of the solution space);
\item If $\rho_p=0$ the term is to be interpreted in the following way:
\begin{equation}\label{dualinterp}
(\langle w_t-\bu \cdot \mathbf e_3,q_p\big|_{s=-h/2}\rangle)_{\omega_p}
=-((w,\partial_t q_p\big|_{s=-h/2}))_{\omega_p}-((\bu \cdot \mathbf e_3,q_p\big|_{s=-h/2}))_{\omega_p}+(w,q_p\big|_{s=-h/2})\big|_{t=0}.
\end{equation}
\end{itemize}

\begin{remark}
We work in the periodic settings to avoid unnecessary technical difficulties. 
However, there are some coercivity issues related to this setting. 
Namely, to ensure coercivity, we require $c_p>0$. 
This is because in the degenerate case $c_p=c_b=0$, the $L^2$ norm of $p_p$ is not bounded by the energy inequality. 
More precisely, in the degenerate case the energy estimate provides the bounds for
 $\|\partial_s p_p\|_{L^2(\Omega_p)}$ and $\|p_p\big|_{s=-h/2}\|_{H^{-1}((0,T)\times \omega_p\times\{-h/2\})}$,
 where the latter bound is a consequence of the coupling condition \eqref{pressure_omega-}, and can be proved by taking $(\bv,0,0,0,0)$ as a test function in 
 the weak formulation \eqref{WeakForm}. 
 However, this is not enough to control high oscillations in the $x_1$ and $x_2$ direction. 
 For example, we can take $p_p=\sin(\pi nx_1)$ to see that its $L^2$ norm is not controlled by the energy.
\end{remark}

{
\begin{remark}
Note that conditions of continuity of the elastic displacement and Biot pressure across the Biot-plate interface are included in the definitions of the solution and test spaces. However, design of finite elements satisfying these constraints can indeed be challenging. One way to address this difficulty is to impose the equality of the traces weakly through Nitsche's method---see for instance, the original paper \cite{nitsche1971variationsprinzip} and \cite{burman2014unfitted} in the context of fluid-structure interactions.
\end{remark}
}

{\begin{remark}\label{testremark} We note that the formal energy inequality in \eqref{Energy} can be obtained by testing with the solution variable in Definition \ref{DefWeakSol} and considering the subsequent cancellations, presented above. \end{remark}}

\subsection{Differential Formulation from Weak Form}
We note that strong solution is recovered from weak solution with sufficient regularity.
Indeed, if a weak solution (as in Definition \ref{DefWeakSol})
$(\bu, w, p_p, \bfeta, p_b)\in \mathcal{V}_{\rm sol}$ is sufficiently smooth in space and time, then 
%
$(\bu, w, p_p, \bfeta, p_b)\in \mathcal{V}_{\rm sol}$ satisfies the PDEs \eqref{Biot_sys}-\eqref{Stokes_sub} in a pointwise sense, as well as the coupling conditions \eqref{kinematic}--\eqref{pressure_omega-}.
Reconstruction of the differential formulation proceeds in the following way:
\begin{enumerate}
\item Take  $\boldsymbol\psi,~q_b,~z,~q_p$  identically zero on their respective domains  in the variational form of Definition 1. 
Invoke the formal integration by parts with $C_c^{\infty}$- test function $\bv$ to observe that the Stokes equation is satisfied point-wise a.e. in $\Omega_f$. 
Then using ``standard'' test functions $\bv$ from the test space ${\cal{V}}_{\rm test}$, the following integral over $\omega_p$ remains:
$$
\int_{\omega_p}\Big[\bfsigma_f(\bu,p_f){\bf e}_3\cdot\bv+\beta(\bu\cdot\bt)(\bv\cdot\bt)
+p_p\big|_{s=-h/2} \bv \cdot \mathbf e_3\Big]=0.
$$
Since $\bv$ is arbitrary, we obtain: 
\begin{enumerate}
\item The Beavers-Joseph-Saffman condition ~$\bfsigma_f{\bf n}\cdot\bt+\beta\bu\cdot\bt=0$ {on $\omega_p$}, by taking $\bv$ to be purely tangential; 
\item The pressure condition $p_p |_{s = -h/2}=-\bfsigma_f\mathbf e_3 \cdot\mathbf e_3$ {on $\omega_p$}, by taking $\bv= \mathbf e_3$.
\end{enumerate}
\item Now take all test function components zero, except {$z$ and $\boldsymbol \psi$}. 
This gives the coupled dynamic equations for the poroelastic plate displacement $w$ and the Biot displacement $\bfeta$:
\begin{align*} 
&-\rho_b((\bfeta_t,\boldsymbol \psi_t))_{\Omega_b}+((\boldsymbol \sigma_b(\bfeta,p_b),\nabla \boldsymbol \psi))_{\Omega_b}
-\rho_p((w_t,\partial_tz))_{\omega_p}+D((\Delta w,\Delta z))_{\omega_p}&
\\[.1cm]\nonumber 
=&\alpha_p((\int_{-h/2}^{h/2} [sp_p]ds,\Delta z))_{\omega_p} 
+ ((p_p\big|_{s=-h/2},z))_{\omega_p}+\rho_b(\bfeta_t,\boldsymbol \psi)\big|_{t=0}+\rho_p(w_t,z)\big|_{t=0}+ ((\mathbf F_b,\boldsymbol \psi))_{\Omega_b}.&
\nonumber
\end{align*}
Similarly as above, to obtain the interior differential equations holding almost everywhere, we use the $C_c^{\infty}$ test functions, and integrate back by parts.
Then, by using the standard test functions $z$ and $\bfpsi$ from ${\cal{V}}_{\rm test}$, after integrating back by parts in the $\bfeta$ equation, we obtain the 
boundary terms on the right hand-side given by 
$$p_p\Big|_{s=-h/2}+\bfsigma_b(\bfeta,p_b)\mathbf e_3\cdot \mathbf e_3.$$

Next, we return to the full weak formulation, including the Stokes equation, and take only $z,~p_p,~\bv$ to be different from zero. 
Undoing the integration by parts in the $\bu$ equation gives the trace terms 
$$-\bfsigma_f(\bu,p_f) \mathbf e_3\cdot \mathbf e_3 +\bfsigma_b(\bfeta,p_b)\mathbf e_3\cdot\mathbf e_3.$$ 
Finally, with the interior equations satisfied a.e., we obtain the following identity holding for the  interface terms, since the appropriate integral equality holds true for all test functions:
$$p_p\Big|_{s=-h/2}+\bfsigma_b(\bfeta,p_b)\mathbf e_3\cdot \mathbf e_3=-\bfsigma_f(\bu,p_f) \mathbf e_3\cdot \mathbf e_3 +\bfsigma_b(\bfeta,p_b)\mathbf e_3\cdot\mathbf e_3.$$
 This is exactly the dynamic coupling condition, i.e. the balance of forces acting across the poroelastic plate.
\item Lastly, let us take $q_p\neq 0$ such that $q_p\big|_{s=h/2}=0$, and $q_b=0$. 
Again, after integrating by parts and taking into account that the interior equations are satisfied point-wise, the following boundary terms remain:
$$
-k_p\int_{\omega_p}\partial_{x_3} p_p\big|_{s=-h/2}q_p\big|_{s=-h/2}+\int_{\omega_p}[w_t-\bu\cdot\mathbf e_3]q_p\big|_{s=-h/2}=0.
$$
Therefore we obtain $\partial_t w=u_3+k_p\partial_3 p_p=u_3-(\bu_p)_3$, where $(\bu_p)_3$ is the third component of the discharge velocity $\bu_p=-k_p\partial_sp_p$,
 and $u_3$ is the third component of the Stokes field $\bu = (u_1,u_2,u_3)$. This is exactly the kinematic coupling condition at the fluid-poroplate interface.
\end{enumerate}
Note that the remaining ``essential" type boundary coupling conditions are built in the compatibility requirements on $\mathcal V_{\rm sol}$ and $\mathcal V_{\rm test}.$

\section{Statement of Main Results and Proof Strategy}
In this section we informally state the following two main results of the paper: 
\begin{enumerate} 
\item The existence of a weak solution to the linear evolution problem: in particular, the Biot coefficient (diagonal tensor) $k_b$ is  a given function of $t$ and $\mathbf x$. 
In this case the following parameters are allowed to be non-negative: $\rho_b,\rho_p, c_b \ge 0$, with only the poroelastic plate's compressibility parameter $c_p>0$.  

\item The existence of a weak solution to the nonlinear quasi-static problem: the coefficient $k_b$ depends on the solution, in particular on the fluid content $\zeta_b=c_bp+\alpha_b\nabla \cdot \eta$.
In this case we consider the  quasi-static Biot model, namely $\rho_b=0$, with the coefficients $\rho_p,c_b \ge 0$ and $c_p>0$.
\end{enumerate}
For the main results, we adopt the semi-discretization approach, also known as Rothe's method, e.g., \cite{zenisek}. In the context of FSI problems, it was developed and used to study moving boundary problems involving linear and nonlinear elastic structures \cite{ARMA,BorisSunNWE}. 
In the present paper, we extend this approach to deal with FPSI problems, showing that this approach is quite robust and applicable to a large class FSI problems. 
Nevertheless, we emphasize that in the present paper,
 the presence of the thin poroelastic plate is of particular importance, since it serves as a regularizing mechanism \cite{reg1,reg2} in the existence proof, while also allowing perfusion/filtration through the interface. 
 
The proof relies on the following main steps:
\begin{enumerate}
\item Construction of approximate solutions.
\begin{enumerate}
\item Derive a semi-discrete in time and linearized version of the weak formulation, with the time discretization parameter $\dt=\frac{T}{N}$. 
In the semi-discretized formulation, the linearization is performed by evaluating the nonlinear coefficient functions at the previous time step. 
This defines a sequence of linear elliptic problems for solutions at time $t^n=n\dt$, with data given at $t^{n-1}=(n-1)\dt$.
\item Prove existence of a unique solution to the elliptic problem defined for every $n=1,\dots,N$ by using the Lax-Millgram Lemma.
\item Define approximate {\emph{solutions}} $$(\bu^{[N]},w^{[N]},p_p^{[N]},\bfeta^{[N]},p_b^{[N]})$$
 as  functions that are piece-wise constant in time, extrapolating the values
obtained from the time discretization, and satisfying the approximate weak formulation.
\end{enumerate}
\item Derive uniform estimates in $\dt$ through appropriate discretized test functions.
\item Obtain weak solution by letting  $N\to\infty$, i.e., $\Delta t \to 0$, utilizing compactness criteria to push limits on the nonlinear term involving $k_b$.
\end{enumerate} 

Since in Steps 1 and 2 we use the linearized problem, these steps are the same for both the linear evolution problem and non-linear quasi-static problem. However, the Step 3 deviates, because in the non-linear case we need to prove additional strong convergence properties. 
In Step 3 we will utilize piecewise constant approximations and an elegant compactness result  for piecewise constant approximate solutions \cite{dreher},
which greatly simplifies the construction from past work on nonlinear Biot problems \cite{zenisek,bgsw,cao}.

We  begin with construction of the approximate solutions that covers both cases. Namely, for the linear problem, simply take $k_b \equiv k_b(\mathbf x,t)$ as a given function.

Finally, we address uniqueness criteria for solutions that exhibit more regularity than our weak solutions. In particular, the result is of weak-strong uniqueness type, wherein weak solutions from a certain class (see Assumption \ref{multiplier}) are unique if there is a single solution in that class exhibiting additional regularity (see Proposition \ref{Uniq}). 

\section{Main Result I: Linear Existence}
In this section we consider $k_b=k_b(\textbf x,t)\in L^{\infty}(0,T;L^{\infty}(\mathbb R))$ as a given function of space and time.

The proof of the existence of a weak solution is based on temporal semi-discretization along with formal energy estimates. 
It is a direct application of Rothe's method to the coupled linear Stokes-Biot-poroelastic plate problem. 
We begin by the construction of approximate solutions, given a discretization of the time interval based on a time step $\Delta t$.
\subsection{Construction of the approximate solutions}

Let $\dt>0=T/N$, $N\in \N$. For every fixed $N$, we inductively define a sequence of approximate solutions $(\bu^n_{N},w^n_{N},(p_p)^n_{N},\bfeta^n_{N},(p_b)^n_{N})$, $n=0,\dots, N$. 
For this purpose, we use implicit Euler method to discretize the time derivatives, making use of the following notation:
 \begin{equation}\label{discderiv}\
 \dot{\bfeta}^{n+1}=\frac{\bfeta^{n+1}-\bfeta^n}{\Delta t},~~~\dot{w}^{n+1}=\frac{w^{n+1}-w^n}{\Delta t}.
 \end{equation}
 Thus, at every fixed time $t^n$, we have to ``solve'' a semi-discretized problem (specified below in \eqref{naiveform}), with the solution (depending only on the spatial variables) 
 belonging to the following {\emph{weak solution space}}:
\begin{align}\label{SDSpace}
\mathcal{V}_{\rm sd}=\Big \{(\bu,w,p_p,\bfeta,p_b) &\in V_{\#}\times H^2_{\#}(\omega_p) \times  H^{0,0,1}\times H^1_{\#}(\Omega_b)  \times H^1_{\#}(\Omega_b)~: \\ &{p_p\big |_{s=h/2}=p_b\big |_{x_3=0},\; \bfeta\big |_{x_3=0}=w{\bf e}_3}\Big \}, \nonumber
\end{align}
equipped with the following norm: 
\begin{align}\label{SDNorm}
\|(\bu,w,p_p,\bfeta,p_b)\|^2_{sd}: &=
\|\bu\|_{H^1(\Omega_f)}^2+\|\Delta w\|_{L^2(\omega_p)}^2+\| w\|_{L^2(\omega_p)}^2
+\|p_p\|^2_{L^2(\Omega_p)}+\|\partial_s p_p\|^2_{L^2(\Omega_p)}\\
&+\|\bfeta\|_{E}^2+\|\nabla p_b\|_{L^2(\Omega_b)}^2,\nonumber
\end{align}
where we recall that the norm $||\cdot||_E$ was defined in \eqref{Enorm}.


\begin{lemma}\label{normequiv}
$(\mathcal{V}_{\rm sd},\|.\|_{sd})$ is a Hilbert space. In particular, the norm $||\cdot||_{\mathcal V_{\rm sd}}$ defined in \eqref{SDNorm}  is equivalent to the usual norm in the space $V_{\#}\times H^2_{\#}(\omega_p) \times  H^{0,0,1}\times H^1_{\#}(\Omega_b)  \times H^1_{\#}(\Omega_b)$.
\end{lemma}

\begin{proof}~
We first notice that $\|\Delta w\|_{L^2(\omega_p)}+\| w\|_{L^2(\omega_p)}$ is a norm equivalent to the $H^2$ norm on $H^2_{\#}(\omega_p)$.
This follows by applying elliptic regularity to the case of periodic boundary conditions associated 
 with the operator $-\Delta w+\gamma w  \in L^2(\Omega)$, giving full control of the $H^2(\omega_p)$ norm.
 
 Next, we observe that the $L^2$ norms of $\bfeta$ and $p_b$ (Biot components) are not included in the definition of $\|.\|_{sd}$,
 and so we need to show the following Poincare-type inequalities:
\begin{align}\label{L2control}
\|\bfeta\|^2_{L^2(\Omega_b)}\leq C\|(\bu,w,p_p,\bfeta,p_b)\|^2_{sd},\quad \bfeta\in H^1_{\#}(\Omega_b),
\\
\|p_b\|^2_{L^2(\Omega_b)}\leq C\|(\bu,w,p_p,\bfeta,p_b)\|^2_{sd},\quad p_b\in H^1_{\#}(\Omega_b).
\end{align}
In contrast with the homogeneous Dirichlet case, the proof of the Poincar\'{e} inequalities will follow by employing the 
interface conditions incorporated in the definition of $\mathcal{V}_{sd}$. 
We will prove the inequalities for smooth functions since general case follows by density argument in the standard way. 
To prove the first inequality in \eqref{L2control} we notice
\begin{align*}
|\bfeta(x_1,x_2,x_3,t)|^2
=|\bfeta(x_1,x_2,0,t)+\int_0^{x_3}\partial_{x_3}\bfeta(x_1,x_2,z,t)dz|^2
\\
\leq C\big (|\bfeta(x_1,x_2,0,t|^2+\int_0^{x_3}|\partial_{x_3}\bfeta(x_1,x_2,z,t)|^2dz \big ).
\end{align*}
Integrate the above inequality over $\Omega_b$ and invoke the interface coupling $\bfeta|_{\omega_p}=\langle 0, 0, w \rangle$ to obtain:
$$
\|\bfeta\|^2_{L^2(\Omega_b)}\leq C\big (\|\bfeta\|_{L^2(w_p)}^2+\|\partial_{x_3}\bfeta\|^2_{L^2(\Omega_b)}\Big )
\leq C\big (\|w\|_{L^2(w_p)}^2+\|\bfeta\|^2_{E}  \big )
\leq C\|(\bu,w,p_p,\bfeta,p_b)\|^2_{sd}.
$$
This proves the first inequality.

To prove the second inequality, we first need to estimate the trace of $p_b$ on $x_3=0$:
$$
p^2_b(x_1,x_2,0,t)=p^2_p(x_1,x_2,h/2,t)=\Big (\frac{2}{h}\int_0^{h/2}\partial_s(sp_p(x_1,x_2,z,t))dz\Big )^2.
$$
We continue, as before, by integrating the above equality over $\omega_p$ to obtain:
$$
\|p_b\|_{L^2(w_p)}^2\leq C\big (\|p_p\|_{L^2(\Omega_p)}^2+\|\partial_s p_p\|_{L^2(\Omega_p)}^2\big).
$$
Thus, we get:
$$
\|p_b\|_{L^2(\Omega_b)}^2
\leq C\big (\|p_b\|_{L^2(w_p)}^2+\|\nabla p_b\|_{L^2(\Omega_b)}^2\big )
\leq C\|(\bu,w,p_p,\bfeta,p_b)\|^2_{sd}.
$$

\end{proof}
 \begin{remark}
If there were homogeneous Dirichlet boundary conditions on any subportion of the boundary for $\bfeta$ and $p_b$, 
the classical Poincar\'e and Korn inequalities would yield norm equivalences for those variables. 
In the present configuration, the $L^2$ control  in \eqref{L2control} for $\bfeta$ and $p_b$ is obtained through interface conditions, 
which gives equivalence of the full norms, as stated in Lemma \ref{normequiv}. 
\end{remark}
 
 We now define weak formulation for the semi-discrete linear problem  at time $t^n=n\dt$ for each $n=1,\dots, N$.
 The weak formulation defines semi-discretized approximations of the solution.
\begin{definition}{\bf{Weak formulation for the semi-discrete problem.}}
The following is the semi-discretized weak formulation of the linear version of problem \eqref{Biot_sys}-\eqref{BiotIC} where $k_b^n=k_b(\mathbf x, t^n)$:
\\
Find $(\bu^n,w^n,(p_p)^n,\bfeta^n,(p_b)^n)\in\mathcal{V}_{\rm sd}$ such that 
for every $(\bv,z,q_p,\bfpsi,q_b)\in \mathcal{V}_{\rm sd}$ the following holds:
\begin{align} 
\nonumber
&\displaystyle{\rho_b(\dot{\bfeta}^{n+1}-\dot{\bfeta}^n,\boldsymbol \psi)_{\Omega_b}
+\dt(\boldsymbol \sigma_b(\bfeta^{n+1},p^{n+1}_b),\nabla \boldsymbol \psi)_{\Omega_b}
+\dt(c_b \dot{p}_b^{n+1}+\alpha_b\nabla \cdot \dot{\bfeta}^{n+1} , q_b)_{\Omega_b}}
&
\\
\nonumber
&
+\dt(k^n_b\nabla p^{n+1}_b,\nabla q_b)_{\Omega_b}
+\dt(c_p\dot{p}_p^{n+1}-\alpha_ps\Delta \dot{w}^{n+1}, q_p)_{\Omega_p}
+\dt(k_p\partial_{s}p^{n+1}_p,\partial_{s}q_p)_{\Omega_p} & 
\\
\nonumber
&+\rho_p(\dot{w}^{n+1}-\dot{w}^n,z)_{\omega_p}
+\dt D(\Delta w^{n+1},\Delta z)_{\omega_p}+\dt\alpha_p(\int_{-h/2}^{h/2} [sp^{n+1}_p]ds,\Delta z)_{\omega_p}&
\\
\nonumber
&+\gamma [\Delta t](w^{n+1},z)_{\omega_p}
-\dt(p^{n+1}_p\Big|_{s=-h/2},z-\bv \cdot \mathbf e_3)_{\omega_p}+\dt(\dot{w}^{n+1}-\bu^{n+1} \cdot \mathbf e_3,q_p\Big|_{s=-h/2})& 
\\
\nonumber
&+ \rho_f(\bu^{n+1}-\bu^n,\bv)_{\Omega_f}+\dt\mu_f (\boldsymbol D(\bu^{n+1}),\boldsymbol D(\bv))_{\Omega_f}+\dt\beta( \bu^{n+1} \cdot \boldsymbol t , \bv \cdot \boldsymbol t  )_{\omega_p}& 
\\ 
\label{naiveform}
=
&~ \dt (\mathbf F^{n+1}_b,\boldsymbol \psi)_{\Omega_b}
+\dt (S^{n+1},q_b)_{\Omega_b}+ \dt (\mathbf f^{n+1},\bv)_{\Omega_f},&
\end{align}
with the initial data specified in \eqref{BiotIC}, and the
discretized time derivative $\dot{\bfeta}^n$ defined in \eqref{discderiv}.
\end{definition}

Notice that in the linear case, $k_b^n$ is the value of $k_b(x,t)$ at the previous time step $t^n$.
In the nonlinear case, we will take $k_b^n=(c_bp_b^n+\alpha_b\nabla\cdot\bfeta^n)$.

To write this weak formulation entirely  in terms of the unknown functions $(\bu^n,w^n,(p_p)^n,\bfeta^n,(p_b)^n)$ belonging to the function space $\mathcal{V}_{\rm sd}$,
we express the terms $\dot{\bfeta}^{n}$ and  $\dot{\bfeta}^{n+1}$ in terms of $\bfeta$ using \eqref{discderiv} to obtain:
\begin{align}
\nonumber
&\rho_b({\bfeta}^{n+1},\boldsymbol \psi)_{\Omega_b}
+{\dt^2}(\boldsymbol \sigma_b(\bfeta^{n+1},p^{n+1}_b),\nabla \boldsymbol \psi)_{\Omega_b}
+\dt(c_b {p}_b^{n+1}, q_b)_{\Omega_b}+\dt(\alpha_b\nabla  {\bfeta}^{n+1} , q_b)_{\Omega_b}
&
\\
\nonumber
&
+\dt^2(k^n_b\nabla p^{n+1}_b,\nabla q_b)_{\Omega_b}
+\dt(c_p{p}_p^{n+1}-\alpha_ps\Delta {w}^{n+1}, q_p)_{\Omega_p}
+\dt^2(k_p\partial_{s}p^{n+1}_p, \partial_{s}q_p)_{\Omega_p} & 
\\
\nonumber
&+\rho_p({w}^{n+1},z)_{\omega_p}
+{\dt^2} D(\Delta w^{n+1},\Delta z)_{\omega_p}+{\dt^2}\alpha_p(\int_{-h/2}^{h/2} [sp^{n+1}_p]ds,\Delta z)_{\omega_p}+\gamma\dt^2(w^{n+1},z)_{\omega_p}&
\\
\nonumber 
&-{\dt^2}(p^{n+1}_p\Big|_{s=-h/2},z)_{\omega_p}+{\dt^2}(p^{n+1}_p\Big|_{s=-h/2} ,  \bv \cdot \mathbf e_3)_{\omega_p}+\dt({w}^{n+1}, q_p\Big|_{s=-h/2})& 
\\
\nonumber
&-\dt^2(\bu^{n+1} \cdot \mathbf e_3, q_p\Big|_{s=-h/2})+ \dt \rho_f(\bu^{n+1},\bv)_{\Omega_f}+\dt^2\mu_f (\boldsymbol D(\bu^{n+1}),\boldsymbol D(\bv))_{\Omega_f}
& 
\end{align}
\begin{align}
 \label{secondform}
&+\dt^2\beta ( \bu^{n+1} \cdot \boldsymbol t , \bv \cdot \boldsymbol t )_{\omega_p}&
\\[.1cm]
\nonumber
=
&~ \dt^2 (\mathbf F^{n+1}_b,\boldsymbol \psi)_{\Omega_b}+\dt^2 (S^{n+1}, q_b)_{\Omega_b}+\dt^2 (\mathbf f^{n+1}, \bv)_{\Omega_f}
\\[.1cm]\nonumber 
&+\dt (\dot{\bfeta}^n,\bfpsi)_{\Omega_b}+\dt (\dot{w}^n,z)_{\omega_p}+\dt \rho_f(\bu^n,\bv)_{\Omega_f}+\rho_b(\bfeta^n,\bfpsi)_{\Omega_b}+\dt(c_b {p}_b^{n}+\alpha_b\nabla  {\bfeta}^{n} , q_b)_{\Omega_b}
\\
\nonumber
&+\dt(c_p{p}_p^{n}-\alpha_ps\Delta {w}^{n}, q_p)_{\Omega_p}+\rho_p({w}^{n},z)_{\omega_p}+\dt({w}^{n}, q_p\Big|_{s=-h/2})
\end{align}
We note that the above problem, although equivalent to the discrete weak form in \eqref{naiveform}, does not have a coercive structure. 
This is a consequence of the particular coupling in our multilayer problem, 
and the use of mixed-velocity/solution test functions in obtaining the energy estimates. 
To circumvent this issue, we modify the test functions (i.e., the variational problem we are trying to solve) in the following way:
\begin{equation}\label{rescaled}
q_b \mapsto {\dt}q_b;~~
q_p \mapsto {\dt}q_p;~~
\bv  \mapsto {\dt}\bv.
\end{equation}
This re-scaling can be undone, of course. 
This is because the re-scaling of test functions does not affect admissability in the space $\mathcal V_{\rm sd}$,
and, after the use of the Lax-Milgram Lemma to obtain existence, as described below, 
the  variational problem in \eqref{secondform} is recovered through the ``for all" quantifier in the definition of a weak solution.
As a consequence, the solution to \eqref{secondform} recovers the solution to the original semi-discretized problem \eqref{naiveform}.

With the re-scaled test functions \eqref{rescaled}, the weak formulation of the semi-discretized linear problem is given by the following:
\begin{definition}\label{weak3}{\bf{Weak formulation with rescaled test functions.}}
We say that a function 
$(\bu^n,w^n,(p_p)^n,\bfeta^n,(p_b)^n) \in \mathcal{V}_{\rm sd}$ is a weak solution to the semi-discretized linear problem
\eqref{Biot_sys}-\eqref{BiotIC} if for every test function $({\dt}\bv,z,{\dt}q_p,\bfpsi,{\dt}q_b)\in \mathcal{V}_{\rm sd}$
the following holds:
\begin{align}
\nonumber
&\rho_b({\bfeta}^{n+1},\boldsymbol \psi)_{\Omega_b}
+ {\dt^2}(\boldsymbol \sigma_b(\bfeta^{n+1},p^{n+1}_b),\nabla \boldsymbol \psi)_{\Omega_b}
+\dt(c_b {p}_b^{n+1},  { \dt}q_b)_{\Omega_b}+\dt(\alpha_b\nabla  {\bfeta}^{n+1} ,  {\dt}q_b)_{\Omega_b}
&
\\
\nonumber
&
+\dt^2(k^n_b\nabla p^{n+1}_b, {\dt}\nabla q_b)_{\Omega_b}
+\dt(c_p{p}_p^{n+1}-\alpha_ps\Delta {w}^{n+1},  {\dt}q_p)_{\Omega_p}
+\dt^2(k_p\partial_{s}p^{n+1}_p, {\dt}\partial_{s}q_p)_{\Omega_p} & 
\\
\nonumber
&+\rho_p({w}^{n+1},z)_{\omega_p}
+ {\dt^2} D(\Delta w^{n+1},\Delta z)_{\omega_p}+ {\dt^2}\alpha_p(\int_{-h/2}^{h/2} [sp^{n+1}_p]ds,\Delta z)_{\omega_p}+\gamma\dt^2(w^{n+1},z)_{\omega_p}&
\\
\nonumber 
&- {\dt^2}(p^{n+1}_p\Big|_{s=-h/2},z)_{\omega_p}+ {\dt^2}(p^{n+1}_p\Big|_{s=-h/2} ,  {\dt} \bv \cdot \mathbf e_3)_{\omega_p}+\dt({w}^{n+1}, {\dt}q_p\Big|_{s=-h/2})& 
\\
\nonumber
&-\dt^2(\bu^{n+1} \cdot \mathbf e_3, {\dt}q_p\Big|_{s=-h/2})& 
\\[.1cm] 
\nonumber
&+  {\dt}\Big[\dt \rho_f(\bu^{n+1},\bv)_{\Omega_f}+\dt^2\mu_f (\boldsymbol D(\bu^{n+1}),\boldsymbol D(\bv))_{\Omega_f}+\dt^2\beta ( \bu^{n+1} \cdot \boldsymbol t , \bv \cdot \boldsymbol t )_{\omega_p}\Big]&
\\[.1cm]
 \label{thirdform}
=
&~ \dt^2 (\mathbf F^{n+1}_b,\boldsymbol \psi)_{\Omega_b}+\dt^2 (S^{n+1}, {\dt}q_b)_{\Omega_b}+\dt^2 (\mathbf f^{n+1}, {\dt}\bv)_{\Omega_f}
\\[.1cm]
\nonumber 
&+\dt (\dot{\bfeta}^n,\bfpsi)_{\Omega_b}+\dt (\dot{w}^n,z)_{\omega_p}+\dt \rho_f(\bu^n, {\dt}\bv)_{\Omega_f}+\rho_b(\bfeta^n,\bfpsi)_{\Omega_b}
+\rho_p({w}^{n},z)_{\omega_p}
\\
\nonumber
&+\dt(c_b {p}_b^{n}+\alpha_b\nabla  {\bfeta}^{n} ,  {\dt}q_b)_{\Omega_b}
+\dt(c_p{p}_p^{n}-\alpha_ps\Delta {w}^{n},  {\dt}q_p)_{\Omega_p}+\dt({w}^{n}, {\dt}q_p\Big|_{s=-h/2}),
\end{align}
with the initial data specified in \eqref{BiotIC}.
\end{definition}

We aim at using the Lax-Milgram lemma to prove existence of a weak solution specified in Definition~\ref{weak3}.
More precisely, we will prove the following result.
\begin{theorem}\label{linear_existence}
There exists a unique weak solution to \eqref{thirdform}.
\end{theorem}
To prove this result using a Lax-Milgram approach we introduce a bilinear form on $\mathcal V_{\rm sd}$, associated with
the left hand-side of \eqref{thirdform}, defining our elliptic operator, and a functional on $\mathcal V_{\rm sd}$ associated with the right hand-side of \eqref{thirdform}
given in terms of $\bfeta^n, \dot\bfeta^n, w^n,  \dot w^n,  \bu^n, p^n_b, p_p^n$ and the source terms $\mathbf F^{n+1}_b,S^{n+1},\mathbf f^{n+1}$, treated as given data.

More precisely, given $k_b^n$, we introduce a bilinear form $a_n(\cdot,\cdot)$ on $\mathcal V_{\rm sd}\times \mathcal V_{\rm sd}$ as:
\begin{align}
\nonumber
a_n&\left((\bfeta, p_b, w, p_p, \bu), (\bfpsi, q_b, z, q_p, \bv)\right) :=
\rho_b({\bfeta}^{n+1},\boldsymbol \psi)_{\Omega_b}
+\dt^2(\boldsymbol \sigma_b(\bfeta^{n+1},p^{n+1}_b),\nabla \boldsymbol \psi)_{\Omega_b}&
\\
\nonumber
&+\dt^2(c_b {p}_b^{n+1}, q_b)_{\Omega_b}+\dt^2(\alpha_b\nabla  {\bfeta}^{n+1} , q_b)_{\Omega_b}
+\dt^3(k_b^n\nabla p^{n+1}_b,\nabla q_b)_{\Omega_b}
&
\\
\nonumber
&
+\dt^2(c_p{p}_p^{n+1}-\alpha_ps\Delta {w}^{n+1}, q_p)_{\Omega_p}
+\dt^3(k_p\partial_{s}p^{n+1}_p,\partial_{s}q_p)_{\Omega_p}  
+\rho_p({w}^{n+1},z)_{\omega_p}
&
\\
\nonumber
&+\dt^2 D(\Delta w^{n+1},\Delta z)_{\omega_p}+\gamma \dt^2 (w^{n+1},z)_{\omega_p},+\dt^2\alpha_p(\int_{-h/2}^{h/2} [sp^{n+1}_p]ds,\Delta z)_{\omega_p}
& 
\\
\nonumber
&
-\dt^2(p^{n+1}_p\Big|_{s=-h/2},z)_{\omega_p}
+\dt^3(p_p^{n+1}\Big|_{s=-h/2},\bv \cdot \mathbf e_3)_{\omega_p}+\dt^2({w}^{n+1},q_p\Big|_{s=-h/2})
& 
\\
\nonumber
&
-\dt^3(\bu^{n+1} \cdot \mathbf e_3,q_p\Big|_{s=-h/2})
+ \dt^2 \rho_f(\bu^{n+1},\bv)_{\Omega_f}+\dt^3\mu_f (\boldsymbol D(\bu^{n+1}),\boldsymbol D(\bv))_{\Omega_f}
&
\\
\label{bilinearform}
&+\dt^3\beta ( \bu^{n+1} \cdot \boldsymbol t , \bv \cdot \boldsymbol t )_{\omega_p}&
\end{align}
Notice again that
the dependence of the bilinear form $a_n(\cdot,\cdot)$ on $n$ 
comes specifically through the linear coefficient $k_b^n$ in the case when this function exhibits any temporal dependence.

Similarly, we introduce the functional $\mathscr F_n$  by:
\begin{align}
\nonumber
\mathscr F_n&\left(\bfpsi, q_b, z, q_p, \bv\right) :=
 \dt^2 (\mathbf F^{n+1}_b,\boldsymbol \psi)_{\Omega_b}+\dt^3 (S^{n+1},q_b)_{\Omega_b}+\dt^3 (\mathbf f^{n+1},\bv)_{\Omega_f}
\\
\nonumber
&+\dt (\dot{\bfeta}^n,\bfpsi)_{\Omega_b}+\dt (\dot{w}^n,z)_{\omega_p}+\dt^2 \rho_f(\bu^n,\bv)_{\Omega_f}+\rho_b(\bfeta^n,\bfpsi)_{\Omega_b}+\dt^2(c_b {p}_b^{n}+\alpha_b\nabla  {\bfeta}^{n} , q_b)_{\Omega_b}
\\
\label{RHSfunx}
&+\dt^2(c_p{p}_p^{n}-\alpha_ps\Delta {w}^{n}, q_p)_{\Omega_p}+\rho_p({w}^{n},z)_{\omega_p}+\dt^2({w}^{n},q_p\Big|_{s=-h/2}),
\end{align}
which depends on $\dt$, the source terms $\mathbf F_b,S,\mathbf f$ at $t^{n+1}$,  and solution evaluated at the previous time step $t^n$.

We have the following lemma:

\begin{lemma}\label{coercivity} Suppose that $c_p>0$ and $\rho_b,\rho_p,c_b\ge 0$.
Consider $a_n: \mathcal V_{\rm sd} \times \mathcal V_{\rm sd} \to \mathbb R$ and $\mathscr F_n \in \mathscr L\big( \mathcal V_{\rm sd}, \mathbb R\big)$, as defined in \eqref{bilinearform} and \eqref{RHSfunx} above. Then, as a bilinear form, $a_n(\cdot,\cdot)$ is continuous and coercive, and $\mathscr F_n \in [\mathcal V_{\rm sd}]'$.
\end{lemma}

Notice that in using Lax-Milgram with our particular rescaling, the setup is robust in that it  allows $c_b, \rho_b, \rho_p$ to vanish. 
\begin{proof}
First we demonstrate coercivity of $a_n(\cdot,\cdot)$ on $\mathcal V_{\rm sd} \times \mathcal V_{\rm sd}$. 
We choose $\bfpsi=\bfeta^{n+1},$ $q_b = p_b^{n+1},$ $q_p =p_p^{n+1}$, and $\bv = \bu^{n+1}$ and calculate:
\begin{align*}
a_n&(\bu^{n+1},w^{n+1},p_p^{n+1},\bfeta^{n+1},p_b^{n+1}),(\bu^{n+1},w^{n+1},p_p^{n+1},\bfeta^{n+1},p_b^{n+1})) = \rho_b||\bfeta^{n+1}||^2_{\Omega_b}+\dt^2 ||\bfeta^{n+1}||^2_E\\
&-\dt^2\alpha_b(p_b^{n+1},\nabla \bfeta^{n+1})_{\Omega_b}+\dt^2(\alpha_b\nabla  {\bfeta}^{n+1} , p^{n+1}_b)_{\Omega_b}+\dt^2 c_b||p_b^{n+1}||^2_{\Omega_b}\\
&+\dt^3||(k^n_b)^{1/2}\nabla p_b^{n+1}||^2_{\Omega_b}+\dt^2c_p||p_p^{n+1}||^2_{\Omega_p}-\dt^2 \alpha_p(s\Delta w^{n+1},p_p^{n+1})_{\Omega_p}+\dt^3||k_p^{1/2}||\partial_sp_p^{n+1}||_{\Omega_p}^2\\
&+\rho_p||w^{n+1}||^2_{\omega_p}+\dt^2 D||\Delta w^{n+1}||^2_{\omega_p} +{\gamma\dt^2 ||w^{n+1}||^2_{\omega_p}}+\dt^2 \alpha_p\left(\int_{-h/2}^{h/2}[sp^{n+1}]ds,\Delta w^{n+1}\right)_{\omega_p}\\
&-\dt^2(p^{n+1}_p\Big|_{s=-h/2},w^{n+1})_{\omega_p}+\dt^3(p_p^{n+1}\Big|_{s=-h/2},\bu^{n+1} \cdot \mathbf e_3)_{\omega_p}+\dt^2({w}^{n+1},p^{n+1}_p\Big|_{s=-h/2}) \\
&-\dt^3(\bu^{n+1} \cdot \mathbf e_3,p^{n+1}_p\Big|_{s=-h/2}) \\\nonumber
&+ \dt^2 \rho_f(\bu^{n+1},\bu^{n+1})_{\Omega_f}+\dt^3\mu_f (\boldsymbol D(\bu^{n+1}),\boldsymbol D(\bu^{n+1}))_{\Omega_f}+\dt^3\beta ( \bu^{n+1} \cdot \boldsymbol t , \bu^{n+1} \cdot \boldsymbol t )_{\omega_p}.
\end{align*}
Now, because of the correct scalings in terms of $\dt$, we can cancel out both Biot type mixed terms, as well as boundary coupling terms, to obtain:
\begin{align*}
a_n&(\bu^{n+1},w^{n+1},p_p^{n+1},\bfeta^{n+1},p_b^{n+1}),(\bu^{n+1},w^{n+1},p_p^{n+1},\bfeta^{n+1},p_b^{n+1}))=\rho_b||\bfeta^{n+1}||^2_{\Omega_b}+\dt^2||\bfeta^{n+1}||^2_E
\\
&+\dt^2 c_b||p_b^{n+1}||^2_{\Omega_b}+\dt^3||[k^n_b]^{1/2}\nabla p_b^{n+1}||^2_{\Omega_b}+\dt^2 c_p||p_p^{n+1}||^2_{\Omega_p}
+\dt^3||[k_p]^{1/2}\partial_sp_p^{n+1}||_{\Omega_p}^2
\\
&+\rho_p||w^{n+1}||^2_{\omega_p}+\dt^2 D||\Delta w^{n+1}||^2_{\omega_p} +{\gamma \dt^2 ||w^{n+1}||^2}
+ \dt^2 ||\bu^{n+1}||^2_{\Omega_f}+\dt^3\mu_f ||\boldsymbol D(\bu^{n+1})||^2_{\Omega_f}
\\
&+\dt^3\beta ||\bu^{n+1} \cdot \boldsymbol t||^2_{\omega_p}.  
\end{align*}
We can now estimate the right hand-side to obtain coercivity on ${\mathcal V_{\rm sd}}$:
\begin{align*}
a_n&(\bu^{n+1},w^{n+1},p_p^{n+1},\bfeta^{n+1},p_b^{n+1}),(\bu^{n+1},w^{n+1},p_p^{n+1},\bfeta^{n+1},p_b^{n+1}))
\\
\gtrsim &~ \dt^2||\bfeta^{n+1}||^2_{\mathbf H^1(\Omega_b)} + \dt^3||p_b^{n+1}||^2_{H^1(\Omega_b)}+\dt^3||p_p^{n+1}||^2_{H^1(\Omega_p)}+\dt^2||w^{n+1}||^2_{H^2(\omega_p)}+\dt^3||\bu||^2_{\mathbf H^1(\Omega_f)}
\\
\gtrsim &~ \text{min}(\dt^2,\dt^3)\|[\bu^{n+1},w^{n+1},p_p^{n+1},\bfeta^{n+1},p_b^{n+1}]\|_{\mathcal V_{\rm sd}}^2.
\end{align*}
Here, we have implicitly used the lower bound $k_b(\mathbf x,t) \ge k_{\rm min}>0$, and the notation 
$\gtrsim$ which represents a bound by a generic constant 
depending on $\rho_f, \rho_b, c_b, \alpha_b, c_p, \alpha_p, \rho_p$, as well as the Korn's and Poincar\'e constants, and possibly $\dt$.

This estimate shows coercivity for the variable $(\bu^{n+1},w^{n+1},p_p^{n+1},\bfeta^{n+1},p_b^{n+1})$ in the $\mathcal{V}_{\rm sd}$ norm.
The quantities $\dot\bfeta^{n+1}$ and $\dot w^{n+1}$ are recovered from the formulae in \eqref{discderiv}, recalling that $\bfeta^n \in \mathbf L^2(\Omega_b)$
and $w^n \in L^2(\omega_p)$ are given as data.

To see that the functional $\mathscr F_n$ belongs to $[\mathcal V_{\rm sd}]'$ we derive a straightforward estimate in terms of $L^2$ norms:
\begin{align*}
||\mathscr F ||_{[\mathcal V_{sd}]'} \lesssim & ~ \dt^2\big[||\mathbf F_b^{n+1}||_{\Omega_b}+\dt ||S^{n+1}||_{\Omega_b} + \dt ||\mathbf f^{n+1}||_{\Omega_f}\big] \\\nonumber
&+\dt\Big[||\dot\bfeta^n||_{\Omega_b}+||\dot w^n||_{\omega_p}+\dt ||\bu^n||_{\Omega_f}\Big]\\\nonumber
&+\dt^2\big[||p^n_b||_{\Omega_b}+||\bfeta^n||_E+||p^n_p||_{\Omega_p}+||\Delta w^n||_{\omega_p}\big] \\\nonumber
& +||\bfeta^n||_{\Omega_b}+||w^n||_{\omega_p}.
\end{align*}
\end{proof}


The proof of Theorem~\ref{linear_existence} now follows by applying the Lax-Milgram lemma to the problem
$$
a_n\left((\bu^{n+1},w^{n+1},p_p^{n+1},\bfeta^{n+1},p_b^{n+1}),(\bfpsi, q_b, z, q_p, \bv)\right) 
= \mathscr F_n \left(\bfpsi, q_b, z, q_p, \bv\right),~~\forall  (\bfpsi, q_b, z, q_p, \bv) \in \mathcal V_{sd},
$$
defined on the Hilbert space $\mathcal V_{sd}$ .

 We emphasize that the approximate solution lies in the test space, as is necessary for obtaining estimates later on.

\begin{remark}
We note that in the above construction, the continuity and coercivity constants may depend on $\dt$ in a singular way. 
However, we utilize the Lax-Milgram lemma only for existence of an approximate solution to the re-scaled problem \eqref{thirdform}. 
After obtaining the solution at the time $t^{n+1}$ from the data at $t^n$, we recover the solution to the original problem \eqref{naiveform} 
Once we have the solution to  \eqref{naiveform}, we obtain {\emph{a priori}} estimates for problem \eqref{naiveform}.
We note that 
the approximate solution lies in the test space, which will be important for the choice of appropriate multiplies in the derivation of {\emph{a priori}} estimates.
The resulting bounds will be uniform in $\dt$. This feature is present in the related considerations \cite{bgsw,zenisek,BW}. 

\end{remark}

Theorem~\ref{linear_existence} provides a sequence of approximate functions, which are functions of the spatial variable only, defined for every $t^n, n = 0,\dots,N$. 
To define  {\emph{approximate solutions}} we extend these functions to the entire time interval $[0,T]$ by extrapolation 
as piecewise constant functions in time, to obtain the following {\emph{approximate solutions}}:
\begin{equation}\label{ApproxSol}
(\bu^{[N]},w^{[N]},p_p^{[N]},\bfeta^{[N]},p_b^{[N]})(t)=(\bu^n_{N},w^n_{N},(p_p)^n_{N},\bfeta^n_{N},(p_b)^n_{N}),\; t\in (t^{n-1},t^n],\; n=1,\dots,N.
\end{equation}
Since we will have to work with derivatives of approximate functions, we introduce the following definition.
\begin{definition}
For a piecewise constant function $f^{[N]}$,  the time derivative of $f^{[N]}$ is a function denoted by $\dot{f}^{[N]}(t)$, defined via backward difference quotients:
\begin{equation}\label{DtN}
\dot{f}^{[N]}(t)=\frac{f^n_N-f^{n-1}_N}{\dt},\quad t\in (t^{n-1},t^n].
\end{equation}
\end{definition}

With this notion of the time derivative, one can easily see that the approximate solutions satisfy the following weak formulation:
\begin{proposition}\label{ApproximateWeakFormulation}
 The approximate solutions $(\bu^{[N]},w^{[N]},p_p^{[N]},\bfeta^{[N]},p_b^{[N]})$, defined by \eqref{ApproxSol}, satisfy the following weak formulation
for all test functions $(\bv,z,q_p,\bfpsi,q_b)\in\mathcal{V}_{test}$:
\begin{align} 
\nonumber
&\rho_b((\ddot{\bfeta}^{[N]},\boldsymbol \psi))_{\Omega_b}
+((\boldsymbol \sigma_b(\bfeta^{[N]},p^{[N]}_b),\nabla \boldsymbol \psi))_{\Omega_b}
+((c_b \dot{p}_b^{[N]}+\alpha_b\nabla \cdot \dot{\bfeta}^{[N]} , q_b))_{\Omega_b}
&\\\nonumber
&
+((k_b(\tau_{\Delta t}(c_bp_b^{[N]}+\alpha_b\nabla\cdot \bfeta^{[N]})\nabla p^{[N]}_b,\nabla q_b))_{\Omega_b}
+((c_p\dot{p}_p^{[N]}-\alpha_ps\Delta \dot{w}^{[N]}, q_p))_{\Omega_p}
+((k_p\partial_{s}p^{[N]}_p,\partial_{s}q_p))_{\Omega_p} & \\\nonumber
&+\rho_p((\ddot{w}^{[N]},z))_{\omega_p}
+ D((\Delta w^{[N]},\Delta z))_{\omega_p}
+\alpha_p((\int_{-h/2}^{h/2} [sp^{[N]}_p]ds,\Delta z))_{\omega_p}&\\\nonumber 
&-((p^{[N]}_p\Big|_{s=-h/2},z-\bv \cdot \mathbf e_3))_{\omega_p}
+((\dot{w}^{[N]}-\bu^{[N]} \cdot \mathbf e_3,q_p\Big|_{s=-h/2}))& \\\nonumber
& + \rho_f((\dot{\bu}^{[N]},\bv))_{\Omega_f}+\mu_f ((\boldsymbol D(\bu^{[N]}),\boldsymbol D(\bv)))_{\Omega_f}
+\beta( \bu^{[N]} \cdot \boldsymbol t , \bv \cdot \boldsymbol t  ))_{\omega_p}& \\ 
\label{semiDform}
&=
 ((\mathbf F^{[N]}_b,\boldsymbol \psi))_{\Omega_b}+ ((S^{[N]},q_b))_{\Omega_b}+ ((\mathbf f^{[N]},\bv))_{\Omega_f},
&
\end{align}
where $\tau_{\Delta t}$ is the 
translation operator defined by $\tau_{\Delta t}f=f(t-\Delta t)$.
\end{proposition}

To see that this is true, we first focus on the time interval $(t^{n},t^{n+1})$ and the weak formulation \eqref{naiveform}.
Let $(\bv,z,q_p,\bfpsi,q_b)\in\mathcal{V}_{test}$, where $\mathcal{V}_{test}$ is defined by \eqref{Testspace}.
Notice that $(\bv,z,q_p,\bfpsi,q_b)(t)$ for $t\in (t^{n},t^{n+1})$,  is an admissible test function for the weak formulation \eqref{naiveform}.
Thus,  take it as a test function in \eqref{naiveform}, integrate the resulting formulation in time from $t^n$ to $t^{n+1}$, 
and then sum with respect to $n=1,\dots,N$, to obtain that the approximate solutions satisfy the approximate weak form \eqref{semiDform}.

\subsection{Uniform {\emph{a priori}} estimates for the approximate solutions}
We show here that the following estimates, uniform in $\dt$, or equivalently, independent of $N$, hold for  approximate solutions.

\begin{lemma}\label{Uniform estimates}
There exists a constant $C > 0$, independent of $N$, such that approximate solutions $\bu^{[N]}$, $w^{[N]}$, $p_p^{[N]}$, $\bfeta^{[N]}$, $p_b^{[N]}$, 
satisfy the following estimates:
\begin{align*}
\rho_f\|\bu^{[N]}\|^2_{L^{\infty}_tL^2_x}
+\rho_p\|\dot{w}^{[N]}\|^2_{L^{\infty}_tL^2_x}
+\rho_b\|\dot{\bfeta}^{[N]}\|^2_{L^{\infty}_tL^2_x}&\leq C,
\\
c_p\|p_p^{[N]}\|^2_{L^{\infty}_tL^2_x}
+\|\bfeta^{[N]}\|_{L^{\infty}_tE}
+c_b\|p_b^{[N]}\|^2_{L^{\infty}_tL^2_x}
+D\|\Delta w^{[N]}\|^2_{L^{\infty}_tL^2_x}&\leq C,
\\
\mu_f\|\boldsymbol D(\bu^{[N]})\|^2_{L^2_tL^2_x}
+\beta\|\bu^{[N]}\cdot {\boldsymbol{t}}\|^2_{L^2_tL^2_x(\omega_p)}
+\|\partial_{x_3}p^{[N]}_p\|^2_{L^2_tL^2_x}
+\|\nabla p^{[N]}_b\|^2_{L^2_tL^2_x}
&\leq C. 
\end{align*}
\end{lemma}
The shorthand notation $L^2_tL^2_x$, for instance, indicates the space $L^2(0,T; L^2(\Omega))$.

\begin{proof}
To simplify notation, in this proof we will be assuming that all the physical constants are equal to 1, except for the constants 
$c_p, c_b, \rho_p, \rho_b$, which can vanish (in the quasistatic case). 

We start by obtaining an {\emph{a priori}} estimate on the functions $\bu^n,w^n,(p_p)^n,\bfeta^n,(p_b)^n$, 
satisfying the semi-discretized problem \eqref{naiveform}.
%
For this purpose, we use the following test functions
in \eqref{naiveform}:  $\bfpsi=\dot\eta^{n+1}$, $q_b=p_b^{n+1}$, $q_p = p_p^{n+1}$, $z=\dot w^{n+1}$ and $\bv= \bu^{n+1}$. 
Using Young's inequality and absorbing all constants that do not depend on $\dt$ we get:
\begin{align*}
&\rho_b||\dot\bfeta^{n+1}||_{\Omega_b}^2+ \dt (\sigma^E(\bfeta^{n+1}),\nabla \dot\bfeta^{n+1})_{\Omega_b}+c_b\dt (\dot p_b^{n+1},p_b^{n+1})_{\Omega_b} +\dt ||(k_b^n)^{1/2}\nabla p_b^{n+1}||_{\Omega_b}^2&\\
&+c_p \dt (\dot p_p^{n+1},p_p^{n+1})_{\Omega_p}+\dt ||k^{1/2}_p\partial_sp_p^{n+1}||^2_{\Omega_p}+\rho_p||\dot w^{n+1}||_{\omega_p}^2+\dt (\Delta w^{n+1},\Delta \dot w^{n+1})_{\omega_p}&\\
&+ ||\bu^{n+1}||^2_{\Omega_f}+\dt  ||\boldsymbol D(\bu^{n+1})||_{\Omega_f}^2 +\dt ||\bu^{n+1}\cdot \boldsymbol t||_{\omega_p}^2\\
\lesssim & ~ \rho_b||\dot \bfeta^n||^2_{\Omega_b} +\rho_p ||\dot w^n||_{\omega_p}^2+||\bu^{n}||_{\Omega_f}^2+ \dt ||\mathbf F^{n+1}_b||^2_{\Omega_b}+\dt ||S^{n+1}||^2_{\Omega_b}+ \dt ||\mathbf f^{n+1}||^2_{\Omega_f}.
\end{align*}
By invoking the definition of difference quotients, as well as Young's inequality, we obtain the following {\emph{a priori}} estimate:
\begin{align*}
&\rho_b||\dot\bfeta^{n+1}||_{\Omega_b}^2+||\bfeta^{n+1}||^2_{E}+c_b||p_b^{n+1}||_{\Omega_b} +\dt ||(k^n_b)^{1/2}\nabla p_b^{n+1}||_{\Omega_b}^2+ c_p||p_p^{n+1}||^2_{\Omega_p}+\dt ||k^{1/2}_p\partial_sp_p^{n+1}||^2_{\Omega_p}& \\
&+\rho_p||\dot w^{n+1}||_{\omega_p}^2+||\Delta w^{n+1}||^2_{\omega_p}+ ||\bu^{n+1}||^2_{\Omega_f}+\dt ||\boldsymbol D(\bu^{n+1})||_{\Omega_f}^2 +\dt ||\bu^{n+1}\cdot \boldsymbol t||_{\omega_p}^2\\
\lesssim & ~  \rho_b||\dot \bfeta^n||^2_{\Omega_b} + ||\bfeta^{n}||_E^2+c_b||p_b^n||_{\Omega_b}^2+c_p||p_p^n||^2_{\Omega_p}+\rho_p ||\dot w^n||_{\omega_p}^2+||\Delta w^n||_{\omega_p}^2 &\\
&+ \dt ||\mathbf F^{n+1}_b||^2_{\Omega_b}+\dt ||S^{n+1}||^2_{\Omega_b}+ \dt ||\mathbf f^{n+1}||^2_{\Omega_f}.&
\end{align*}
This estimate can be rearranged and expressed as follows:
\begin{align*}
&\rho_b||\dot\bfeta^{n+1}||_{\Omega_b}^2+||\bfeta^{n+1}||^2_{E}+c_b||p_b^{n+1}||_{\Omega_b} + c_p||p_p^{n+1}||^2_{\Omega_p}+ \rho_p ||\dot w^{n+1}||_{\omega_p}^2+ ||\Delta w^{n+1}||^2_{\omega_p}+ ||\bu^{n+1}||^2_{\Omega_f}& \\
\lesssim & ~  \rho_b||\dot \bfeta^n||^2_{\Omega_b} + ||\bfeta^{n}||_E^2+c_b||p_b^n||_{\Omega_b}^2+c_p||p_p^n||^2_{\Omega_p}+ \rho_p||\dot w^n||_{\omega_p}^2+||\Delta w^n||_{\omega_p}^2&\\ &+ \dt ||\mathbf F^{n+1}_b||^2_{\Omega_b}+\dt ||S^{n+1}||^2_{\Omega_b}+ \dt ||\mathbf f^{n+1}||^2_{\Omega_f} &\\
&-\dt ||k_b^{1/2}\nabla p_b^{n+1}||_{\Omega_b}^2-\dt ||k^{1/2}_p\partial_sp_p^{n+1}||^2_{\Omega_p}-\dt ||\boldsymbol D(\bu^{n+1})||_{\Omega_f}^2-\dt ||\bu^{n+1}\cdot \boldsymbol t||_{\omega_p}^2.&
\end{align*}
Inductively, after rearrangement, we obtain:
\begin{align}\nonumber
||&\dot\bfeta^{N+1}||_{\Omega_b}^2+||\bfeta^{N+1}||^2_{E}+c_b||p_b^{N+1}||_{\Omega_b} + c_p||p_p^{N+1}||^2_{\Omega_p}+\rho_p||\dot w^{N+1}||_{\omega_p}^2+||\Delta w^{N+1}||^2_{\omega_p}+ ||\bu^{N+1}||^2_{\Omega_f}&\\ \nonumber
&+\sum_{n=0}^{N}\Big[||k_b^{1/2}\nabla p_b^{n+1}||_{\Omega_b}^2+ ||k^{1/2}_p\partial_sp_p^{n+1}||^2_{\Omega_p}+ ||\boldsymbol D(\bu^{n+1})||_{\Omega_f}^2+||\bu^{n+1}\cdot \boldsymbol t||_{\omega_p}^2 \Big] \dt & \\\nonumber
\lesssim &~ ||\dot \bfeta^0||^2_{\Omega_b} + ||\bfeta^{0}||_E^2+c_b||p_b^0||_{\Omega_b}^2+c_p||p_p^0||^2_{\Omega_p}+ ||\dot w^0||_{\omega_p}^2+\rho_p||\Delta w^0||_{\omega_p}^2&\\ &+ \sum_{n=0}^{N} \Big[||\mathbf F^{n+1}_b||^2_{\Omega_b}+||S^{n+1}||^2_{\Omega_b}+ ||\mathbf f^{n+1}||^2_{\Omega_f} \Big]\dt.&
\end{align}
Noting the lower bound on $k_b \ge k_{\rm min}>0$ again,  we obtain the uniform estimates for the approximate solutions.

\end{proof}

We will need two more estimates to be able to pass to the limit as $N \to \infty$ (or, equivalently, as $\Delta t \to 0$) 
in the evolution case when we have to deal with the terms $\dot \bfeta^{[N]}$ and $\dot w^{[N]}$.
In this case we  require bounds on $\ddot \bfeta^{[N]}, ~\ddot w^{[N]}$ in the associated dual, $H^{-1}$-like spaces. 
Additionally, we will need an estimate for the time derivative of the fluid content in the associated dual space.

Having established the solution, such estimates will follow immediately from \eqref{naiveform}. 
More precisely, since some test functions are coupled, we need to consider duals of  product spaces to establish these estimates.
In particular, consider first taking all test functions in \eqref{naiveform} zero except for the ``inertial test functions''
$(\bfpsi,z) \in  H^1_{\#}(\Omega_b)\times H^2_{\#}(\omega_p)$ such that $\bfpsi|_{x_3=0}=z{\bf e}_3$. 
Then, dividing \eqref{naiveform} by $\dt$, we obtain the following equality for the solution at any $t^{n+1}$:
\begin{align*}  
&\rho_b(\ddot \bfeta^{n+1},\bfpsi)_{\Omega_b}
+\rho_p (\ddot w^{n+1},z)_{\omega_p}=  {{(\boldsymbol \sigma_b(\bfeta^{n+1},p_b^{n+1}),\nabla \bfpsi))_{\Omega_b}}}+
 (\mathbf F^{n+1}_b,\boldsymbol \psi)_{\Omega_b}
 \\
 &- D(\Delta w^{n+1},\Delta z)_{\omega_p}-\alpha_p(\int_{-h/2}^{h/2} [sp^{n+1}_p]ds,\Delta z)_{\omega_p}-(p^{n+1}_p\Big|_{s=-h/2},z)_{\omega_p}.
\end{align*}
Similar calculation can be done for the fluid content by taking all test functions zero except for the ``fluid content test functions''
$(q_p,q_b)\in H^1_{\#}(\Omega_b)\times H^{0,0,1}_{\#}$ such that $q_b|_{s=h/2}=q_p$. 
After integrating the obtained equalities with respect to time from $t^n = n \dt$ to $t^{n+1} = {(n+1)\dt}$, 
and taking the sum over $n$, for $n=0,\dots, N-1$,  the estimates in Lemma \ref{Uniform estimates} give the following result:

\begin{lemma}\label{UniformDualEstimates}
Let $\mathcal{V}_{\rm in}=\{(\bfpsi,z) \in  H^1_{\#}(\Omega_b)\times H^2_{\#}(\omega_p):\bfpsi|_{x_3=0}=z{\bf e}_3\}$ and 
$\mathcal{V}_{\rm fc}=\{(q_p,q_b)\in H^1_{\#}(\Omega_b)\times H^{0,0,1}:q_b|_{s=h/2}=q_p\}$. Then there exists constant $C$ independent of $N$ such that:
\begin{align*}
\Big\|\big (\rho_b\ddot{\bfeta}^{[N]},\rho_p\ddot{w}^{[N]}\big )\Big\|_{\mathcal{V}_{\rm in}'}\leq C,
\\
\Big\|\big (c_b\dot{p}_b^{[N]}+\alpha_b\nabla\cdot\dot{\bfeta}^{[N]},c_p\dot{p}_p^{[N]}-\alpha_p\Delta\dot{w}^{[N]}\big )\Big \|_{\mathcal{V}_{\rm fc}'}\leq C.
\end{align*}
\end{lemma}

\subsection{Limit $N\to\infty$ and Recovery of Linear Weak Solution}
\subsubsection{Existence} To prove the existence of a weak solution to the linear  Stokes-Biot-poroplate fluid-structure interaction problem \eqref{Biot_sys}--\eqref{BiotIC},
we would like to show that there exist subsequences of approximate solutions $\{(\bu^{[N]},w^{[N]},p_p^{[N]},\bfeta^{[N]},p_b^{[N]})\}_{N\in\N}$, which
converge,  as $N \to \infty$ (or equivalently, as $\Delta t \to 0$), to a function which satisfies the weak formulation of the continuous problem
specified in Definition \ref{DefWeakSol}. 
Since the problem is linear, weak convergence of subsequences of approximate solutions 
$\{(\bu^{[N]},w^{[N]},p_p^{[N]},\bfeta^{[N]},p_b^{[N]})\}_{N\in\N}$ is sufficient to pass to the limit in the weak formulations \eqref{semiDform} satisfied
by the approximate functions, 
and recover the  weak formulation of the continuous problem, satisfied by the limiting function.
Indeed, the uniform bounds presented in Lemmas \ref{Uniform estimates} and \ref{UniformDualEstimates} imply existence of subsequences
of $\{(\bu^{[N]},w^{[N]},p_p^{[N]},\bfeta^{[N]},p_b^{[N]})\}_{N\in\N}$, denoted again with $\{(\bu^{[N]},w^{[N]},p_p^{[N]},\bfeta^{[N]},p_b^{[N]})\}_{N\in\N}$
(with a slight abuse of notation), such that the following weak convergence is true for those subsequences:
\begin{equation}\label{WeakLimits}
(\bu^{[N]},w^{[N]},p_p^{[N]},\bfeta^{[N]},p_b^{[N]})\rightharpoonup (\bu,w,p_p,\bfeta,p_b)\quad{\rm weakly\; in}\; \mathcal{V}_{sol},
\end{equation}
and
\begin{align*}
(\rho_b\ddot{\bfeta}^{[N]},\rho_p\ddot{w}^{[N]}\big )& \overset{*}\rightharpoonup (\rho_b\ddot{\bfeta},\rho_p\ddot{w}\big )\quad{\rm weakly*\; in}\; \mathcal{V}_{\rm in}
\\
(c_b\dot{p}_b^{[N]}+\alpha_b\nabla\cdot\dot{\bfeta}^{[N]},c_p\dot{p}_p^{[N]}-\alpha_p\Delta\dot{w}^{[N]}\big ),
& \overset{*}\rightharpoonup 
(c_b\dot{p}_b^+\alpha_b\nabla\cdot\dot{\bfeta},c_p\dot{p}_p-\alpha_p\Delta\dot{w}\big )\quad{\rm weakly*\; in}\; \mathcal{V}_{\rm fc}.
\end{align*}
This is sufficient to pass to the limit in \eqref{semiDform}, 
and obtain that the limiting function satisfies the continuous weak formulation \eqref{WeakForm}.
This calculation is straight forward. 
Therefore, the uniform estimates from Lemmas \ref{Uniform estimates} and \ref{UniformDualEstimates} imply the following existence result:

\begin{theorem}\label{MainLinear}
Let $c_p>0$ and $\rho_b,\rho_p,c_b,\rho_f\geq 0$, $\alpha_n,\alpha_p,\mu_E,\lambda_E,D,\mu,k_p>0$, and $T>0$. Moreover,  let 
the initial data belong to the following spaces: {
$
{\bf u} \in L^2(\Omega_f), 
w \in H_{\#}^2(\omega_p), w_t \in L^2(\omega_p),
{\boldsymbol\eta} \in H_{\#}^1(\omega_p), \eta_t \in L^2(\Omega_b),
p_b \in L^2(\Omega_b),
p_p \in L^2(\Omega_p),
$}
and let $k_b\in L^{\infty}(0,T;L^{\infty}(\Omega_b))$ such that $0<k_{\rm min}\leq k_b\leq k_{\rm max}$. Then there exists a weak solution on $(0,T)$ in the sense of Definition \ref{DefWeakSol}.
\end{theorem}

\subsubsection{Uniqueness} We conclude this section with a brief discussion of the issue of (linear) uniqueness. First, we note that even in the linear case, uniqueness of solutions is a non-trivial issue---see \cite{BW,showmono}. This primarily stems from the fact that the energy estimates obtained for the constructed solution above need not apply to a general weak solution. Hence, given two weak solutions, an {\emph{a priori}} estimate on the difference is not immediately obtainable via formal estimation as was done on 
approximate solutions in the construction of our particular weak solution. This issue is two-pronged: first, the hyperbolic nature of the inertial problem (when $\rho_b$ and/or $\rho_p$ non-zero) has challenges associated to justifying velocity multipliers (test functions) $\bfeta_t$ and $w_t$, owing to the higher order terms appearing from the Biot structure of the solid dynamics. Secondly, Biot models with time dependent coefficients (even in the quasi-static case) have issues associated to multiplier analysis in obtaining a priori estimates. The effect of this can be noted, for instance, in \cite[pp.110--117]{showmono} where abstract implicit evolutions are discussed and ``much non-uniqueness" is possible, since ``...the stability estimate...was shown to hold for {\em some} solution, not every solution...."

We provide a lemma that yields uniqueness, assuming that formal test functions $\bfeta_t$ and $w_t$ are admissible in the sense of the following hypothesis:
\begin{assumption}\label{multiplier}
Assume that the function $(\bu,w,p_p,\bfeta,p_b) \in \mathcal V_{\rm sol}$ has the additional property that 
\begin{align}\label{mult1}
( \rho_b w_{tt}+D\Delta^2 w,  w_t)_{\omega_p}  =&~\frac{1}{2}\dfrac{d}{dt}\Big\{\rho_p||w_t||^2_{L^2(\omega_p)}+D||\Delta w||^2_{L^2(\omega_p)} \Big\}\\ \label{mult2}
( \rho_p\bfeta_{tt}+\boldsymbol \sigma^E, \bfeta_t)_{\Omega_p} =&~  \frac{1}{2}\dfrac{d}{dt}\Big\{\rho_b||\bfeta_t||^2_{L^2(\Omega_b)}+||\bfeta||_E^2\Big\}.
\end{align}
\end{assumption}
Then, the following uniqueness result holds.
\begin{lemma}[Uniqueness of Weak Solutions]\label{linuniq}
Let $c_p>0$ and $\rho_b,\rho_p,c_b,\rho_f\geq 0$, $\alpha_b,\alpha_p,\mu_E,\lambda_E,D,\mu,k_p>0$, and $T>0$. Moreover,  let $k_b\in L^{\infty}(\Omega_b\times (0,T))$ such that $0<k_1\leq k_b\leq k_2$. Suppose further that a weak solution $(\bu,w,p_p,\bfeta,p_b) \in \mathcal V_{\rm sol}$ has the additional property specified in Assumption \ref{multiplier}. 
Then, for the given data, such a solution is unique.
\end{lemma}
\begin{proof}
The above lemma follows immediately from the formal energy estimates presented in the preceding section applied to the difference of two solutions to the linear problem. 
\end{proof} 
{\begin{remark} 
Note that Assumption \ref{multiplier} holds, for instance, for strong solutions (for instance if $w \in C^2([0,T];L^2(\omega_p))\cap C^1([0,T];H^2(\omega_p)) \cap C([0,T];H^4(\omega_p))$, with analogous assumptions for $\bfeta$). However, it is not at all obvious when, for instance, $\rho_p$ or $\rho_b$ are zero. 
In full generality, whether or not $w_t$ and $\bfeta_t$ are valid to test with for a weak solution is certainly a subtle issue.\end{remark}}

%
%
%
%
%
%

\section{Main Result II: Quasistatic Nonlinear Problem}\label{sec:quasistatic}

In this section we consider the nonlinear case where:
\begin{itemize}
 \item the permeability function $k_b$ depends explicitly on 
$\zeta_b=c_b p_b + \alpha_b \nabla \cdot \boldsymbol\eta$, as specified in \eqref{nonperm};
\item the Biot model \eqref{Biot_sys} is quasi-static, 
as specified in Assumption~\ref{quasistatic}, namely, 
$\rho_b = 0$ in \eqref{Biot_sys};
\item  we allow $\rho_p\geq 0$ in the poroelastic plate equations \eqref{plate_sys}.
\end{itemize}

The main nonlinear existence result is the following:
\begin{theorem}[Nonlinear existence of weak solutions]\label{MainNonLin}
Let $\rho_b=0$, $c_p>0$, and $\rho_p,c_b,\rho_f\geq 0$, $\alpha_n,\alpha_p,\mu_E,\lambda_E,D,\mu,k_p>0$, and $T>0$. Moreover, let the initial data belong to the following spaces:
$
{\bf u} \in L^2(\Omega_f), 
w \in H_{\#}^2(\omega_p), w_t \in L^2,
{\boldsymbol\eta} \in H_{\#}^1(\Omega_b), \eta_t \in L^2(\Omega_b),
p_b \in L^2(\Omega_b),
p_p \in L^2(\Omega_p),
$
and let $k_b(\cdot) \in L^{\infty}(\mathbb R)$ such that $0<k_{\rm min}\leq k_b\leq k_{\rm max}$. Then, there exists a weak solution to the nonlinear quasi-static problem on $(0,T)$ in the sense of Definition \ref{DefWeakSol}, with $\rho_b=0$.
\end{theorem}

To prove this result, we first notice that all the steps in the existence proof performed for the linear case still hold, except the last step involving passage to the limit. 
To pass to the limit in the nonlinear problem we will use a version of Aubin-Lions lemma obtained in \cite{dreher},
which is particularly suitable for piece-wise constant approximations in time of the approximate solutions
constructed using Rothe's method. 
Here we state the special case of \cite[Theorem 1]{dreher} adapted to our notation:

\begin{theorem}[Compactness Lemma \cite{dreher}]\label{DreherCompact}
Let $V \hookrightarrow H \hookrightarrow W$ be Hilbert spaces with the first inclusion compact and the second continuous. 
Consider the sequence of piece-wise constant in time functions $u^{[N]}$ (indexed by $N\to \infty$, or equivalently, $[\Delta t] \to 0$),
and $\dot{u}^{[N]}$ defined in \eqref{DtN}.

If the following inequality holds
$$\big|\big|\dot{u}^{[N]}\big|\big|_{L^2([\Delta t],T;W)} + ||u^{[N]}||_{L^2(0,T;V)} \le C,$$
 with $C$ independent of $N$, then $u^{[N]}$ is pre-compact in $L^2(0,T;H)$.
\end{theorem}

We will apply this theorem on the sequence $\big[c_bp_b^{[N]}+\alpha_b\nabla \cdot \bfeta^{[N]}\big]$,
with $W = L^2\big(0,T;H^{-1}(\Omega_b)\big)$, $V = L^2\big(0,T;H^{\epsilon}(\Omega_b)\big)$, and $H = L^2\big(0,T;L^2(\Omega_b)\big)$.


We start by showing the following uniform estimate on $\big[c_b\dot{p_b}^{[N]}+\alpha_b\nabla\cdot\dot{\bfeta}^{[N]}\big]$:
\begin{lemma}\label{FCDual}
There exists a constant $C$ independent of $N$ such that
$$
\|c_b\dot{p_b}^{[N]}+\alpha_b\nabla\cdot\dot{\bfeta}^{[N]}\|_{L^2(0,T;H^{-1}(\Omega_b))}\leq C.
$$
\end{lemma}
\proof
Let $q_b\in C^1_c((0,T);H^1_0(\Omega_b))$. Then $(0,0,0,0,q_b)$ belongs to the space of test functions ${V}_{test}$. 
Therefore we can use this test function in the weak formulation \eqref{semiDform} for approximate solutions  to obtain
(due to the second equation in \eqref{Biot_sys}):
$$
((c_b \dot{p}_b^{[N]}+\alpha_b\nabla \cdot \dot{\bfeta}^{[N]} , q_b))_{\Omega_b}
=-((k_b(\tau_{\Delta t}(c_bp_b^{[N]}+\alpha_b\nabla\cdot \bfeta^{[N]}))\nabla p^{[N]}_b,\nabla q_b))_{\Omega_b}
+((S^{[N]},q_b))_{\Omega_b}.
$$
By using uniform estimates presented in Lemma \ref{Uniform estimates}, we obtain the following estimate:
$$
|((c_b \dot{p}_b^{[N]}+\alpha_b\nabla \cdot \dot{\bfeta}^{[N]} , q_b))_{\Omega_b}|
\leq k_{max}\|\nabla p_b^{[N]}\|_{L^2_tL^2_x}\|q_b\|_{L^2_tH^1_x}+C\|S^{[N]}\|_{L^2_tL^2_x}\|q_b\|_{L^2_tH^1_x}
\leq C\|q_b\|_{L^2_tH^1_x}.
$$
Now the conclusion of the lemma follows by density arguments and from the converse of H\"older's inequality.
\qed

In order to prove the compactness result, we need the regularity result for the Biot displacement.
\begin{theorem}\label{RegTm}
In the quasistatic case ($\rho_b=0$), the Biot displacement satisfies the following additional regularity property:
the sequence $\bfeta^{[N]}$ is uniformly bounded in $L^2(0,T;\mathbf H^{2}(\Omega_b))$.
\end{theorem}
\proof
We use the arguments and notation similar to \cite[Section 3.3]{bgsw}. 
We denote by $\mathcal{E}: \mathcal D(\mathcal E) \subset \mathbf L^2(\Omega_b) \to \mathbf L^2(\Omega_b)$ 
the Lam\'{e} operator:
$$
\mathcal{E}(\bfeta) := - \nabla \cdot \left( 
2\mu_b \boldsymbol{D}(\boldsymbol\eta) + \lambda_b \nabla \cdot \boldsymbol\eta \boldsymbol{I}
\right)
$$
defined on the domain
$$
\mathcal D(\mathcal E) = \{
\boldsymbol\eta \in {\bf{H}}^1_{\#}(\Omega_b) ~|~ \nabla \cdot \left( 
2\mu_b \boldsymbol{D}(\boldsymbol\eta) + \lambda_b \nabla \cdot \boldsymbol\eta \boldsymbol{I}
\right)
\in \mathbf{L}^2(\Omega_b)
\}.
$$

Then $\bfeta^{[N]}$ satisfies the following equation:
$$
\mathcal{E}\bfeta^{[N]}=-\nabla p^{[N]}\quad {\rm in}\quad (0,T)\times \Omega_b,
$$
with the following mixed boundary conditions:
\begin{align*}
&{\rm 1.\ Periodic} \ {\rm at}\ x_i = 0, x_i = 1, \ {\rm for} \ i = 1,2 \ {\rm(see\  \eqref{PeriodicBC})};\\
&{\rm 2.\ Zero\  traction\ condition \ at} \ x_3 = 1 \ {\rm(see\  \eqref{zero_traction})};\\
& {\rm 3.\ Dirichlet \ condition \ on} \ x_3 = 0, \ {\rm i.e. \ on\ } \omega_p:\  \bfeta^{[N]}=w^{[N]}{\bf e}_3 \ {\rm{(see \ \eqref{kinematic})}}.
\end{align*}
From Lemma \ref{Uniform estimates} we see that:
\begin{itemize}
\item  The right hand-side $\nabla p$ in uniformly bounded in $L^2(0,T;L^2(\Omega_b))$;
\item The sequence $w^{[N]}$ in the Dirichlet boundary condition is uniformly bounded in $L^{\infty}(0,T;H^2(\omega_p))$.
\end{itemize}
Now the conclusion of the Lemma follows directly from standard elliptic regularity on a rectangular prism, as in \cite[Lemma 2]{bgsw}.
\qed
\begin{remark}
First, it is important to note here that due to the way how the problem is set up, our domain is smooth without corners. Secondly, the full elliptic regularity recovered in Theorem \ref{RegTm} is more than is needed in the construction of solutions; indeed, if $\bfeta^{[N]}$ is uniformly bounded in $L^2(0,T;\mathbf H^{1+\epsilon}(\Omega_b)),~~\epsilon>0,$ then the 
convergence properties stated below follow as in \cite{bgsw}.
\end{remark}

\begin{corollary}\label{StongConvergence}
The following convergence properties hold:
\begin{align*}
c_b{p_b}^{[N]}+\alpha_b\nabla\cdot{\bfeta}^{[N]}\to c_b p_b+\alpha_b\nabla\cdot\bfeta\quad{\rm strongly\;\rm in}\; L^2(0,T;L^2(\Omega_b))
\\
k_b\big (\tau_{\Delta t}(c_bp_b^{[N]}+\alpha_b\nabla\cdot \bfeta^{[N]}) \big )\to k_b\big (c_b p_b+\alpha_b\nabla\cdot\bfeta\big )\quad{\rm strongly\;\rm in}\; {L^2(0,T;L^2(\Omega_b))}
\end{align*}
\end{corollary}
\proof
Uniform energy estimates from Lemma \ref{Uniform estimates} imply that sequence $p_b^{[N]}$ is uniformly bounded in $L^2(0,T;H^1(\Omega_b))$. Therefore, by Theorem \ref{RegTm}, sequence $c_b{p_b}^{[N]}+\alpha_b\nabla\cdot{\bfeta}^{[N]}$ is uniformly bounded in $L^2(0,T;H^{\epsilon}(\Omega_b))$. 
Moreover, by Lemma \ref{FCDual} sequence 
$\Big (c_b{ p_b}^{[N]}+\alpha_b\nabla\cdot{\bfeta}^{[N]}\Big )^{\cdot}$ is uniformly bounded in $L^2(0,T;H^{-1}(\Omega_b))$. 

The strong convergence of $\big[ c_b{p_b}^{[N]}+\alpha_b\nabla\cdot{\bfeta}^{[N]}\big]$ then follows by direct application of Theorem \ref{DreherCompact}.

The second convergence result follows by noticing that $k_b(\cdot): L^2(\Omega_b) \to L^2(\Omega_b)$ represents a Nemytskii operator \cite{showmono}.
This follows from the bounds $0<k_{\rm min} \le k_b(s) \le \kappa_{\rm max}$ for all $s \in \mathbb R$, as well as the continuity of $k_b(\cdot)$. Thus 
$$k_b\big (\tau_{\Delta t}(c_bp_b^{[N]}+\alpha_b\nabla\cdot \bfeta^{[N]}) \big )\to k_b\big (c_b p_b+\alpha_b\nabla\cdot\bfeta\big ),$$ as desired.
\qed

To complete the proof of the existence of weak solutions, what remains to show is how to pass to the limit in the nonlinear term in \eqref{semiDform}. 
For this purpose, we recall the following convergence properties of the sequences appearing in the nonlinear term:
$$
k_b\big (\tau_{\Delta t}(c_bp_b^{[N]}+\alpha_b\nabla\cdot \bfeta^{[N]}) \big )\to k_b\big (c_b p_b+\alpha_b\nabla\cdot\bfeta\big )\;{\rm in}\; L^2(0,T;L^2(\Omega_b)),
$$
$$
\nabla p_b^{[N]}\rightharpoonup \nabla p_b\quad{\rm weakly\;\rm in}\; L^2(0,T;L^2(\Omega_b)).
$$
Now, by the convergence of weak-strong products we get: 
$$
k_b\big (\tau_{\Delta t}(c_bp_b^{[N]}+\alpha_b\nabla\cdot \bfeta^{[N]} )\big )\nabla p_b^{[N]}\rightharpoonup 
 k_b\big (c_b p_b+\alpha_b\nabla\cdot\bfeta\big ) \nabla p_b~~\text{in}~~L^1(0,T;L^1(\Omega_b)).
$$
To see that this sequence converges weakly in $L^2(0,T;L^2(\Omega_b))$, we recall that  $k_b\big (\tau_{\Delta t}(c_bp_b^{[N]}+\alpha_b\nabla\cdot \bfeta^{[N]}) \big )\nabla p_b^{[N]}$ 
is a bounded sequence in $L^2(0,T;L^2(\Omega_b))$, and therefore it has a subsequence which converges weakly in $L^2(0,T;L^2(\Omega_b))$. 
Because of the uniqueness of the limit in $\mathcal{D}'$, its limit is exactly $k_b\big (c_b p_b+\alpha_b\nabla\cdot\bfeta\big ) \nabla p_b$.

\qed
\subsection{Uniqueness}
Lastly, we  consider the issue of uniqueness for nonlinear solutions. At present, even for the nonlinear Biot system alone, 
a full uniqueness result does not exist in the literature. 
In \cite{cao}, a {\em central} hypothesis for the entire well-posedness analysis is that $c_b>0$;  for uniqueness, the permeability function $k_b$ is assumed to be Lipschitz continuous, {\em and} a smallness hypothesis is also imposed on the Biot pressure.

Along these lines, we present here a strong-weak uniqueness result  for mild solutions. More precisely, 
in addition to assuming that our weak solutions satisfy an energy-type equality, specified in Assumption \ref{multiplier}, which is related 
to allowing $\boldsymbol\eta_t$ to serve as a test function (multiplier), we also require some additional regularity on the Biot pressure, as specified below in Lemma~\ref{Uniq}.

More precisely, we introduce the space
$$\mathcal W \subset \mathcal V_{\rm sol} \ {\rm such \ that \  Assumption \ \ref{multiplier} \ is \ satisfied,}$$
and consider solutions in $\mathcal W$ with the additional regularity for $p_b$ as follows:

\begin{proposition}[Uniqueness of Weak Solutions]\label{Uniq}
Let $c_p>0$ with $\rho_b=0$, and $\rho_p,c_b,\rho_f\geq 0$, $\alpha_n,\alpha_p,\mu_E,\lambda_E,D,\mu,k_p>0$. Moreover,  let $k_b\in L^{\infty}(\mathbb R)\cap Lip(\mathbb R)$ such that $0<k_{\rm min}\leq k_b\leq k_{\rm max}$. Suppose there is a weak solution $(\bu,w,p_p,\bfeta,p_b) \in \mathcal W$ and a time $T^*>0$ so that the following holds:
$$ 
\nabla{{p_b}} \in L^2(0,T^*;\mathbf{L}^{\infty}(\Omega_b)).
$$
\noindent
Then, all weak solutions in $\mathcal W$ are equal to $(\bu,w,p_p,\bfeta,p_b)$ on any interval $t\in [0,T]$ with $T<T^*$.
\end{proposition}

\begin{proof} 
Assume that there are two solutions
$\mathbf s^i=(\bu^i,w^i,p_p^i,\bfeta^i, p_b^i) \in \mathcal W$ for $i=1,2$.
Under the  hypotheses {{that solutions belong to $\mathcal W$}}, each solution satisfies the 
formal energy inequality \eqref{Energy}.
Consider the difference $\overline{\mathbf s}=\mathbf s^1-\mathbf s^2$ (with the ``bar" notation for each coordinate as well) 
as a weak solution to the {\emph{difference}} equation with the same data $\mathbf F_b, S, F_p , \mathbf f$, utilizing (in the appropriate sense) $\overline{\mathbf s}$ as a test function. 
Then, for any $T < T^*$, $\overline{\mathbf s}$ satisfies the inequality: 
\begin{align}\nonumber
\mathcal E(\overline{\mathbf s}(T)) +&{{\int_0^T\Big[
{2}\mu_f\|\boldsymbol{D}(\overline{\bu}(t))\|^2_{L^2(\Omega_f)}+ \beta\|\overline{\bu}(t)\cdot\bt\|^2_{L^2({\omega_p})}+
\|k_p^{1/2}\partial_{s} \overline{p}_p(t)\|^2_{L^2(\Omega_p)}\Big]}} dt \\
\le &~  \mathcal E(\overline{\mathbf s}(0))-\int_0^T \left([k_b^1(t)\nabla p_b^1(t)-k_b^2(t)\nabla p_b^2(t)],\nabla [p_b^1(t)-p_b^2(t)]\right)_{\Omega_b}dt
\end{align}
where we denoted
$$
\mathcal E(\mathbf s^i) \equiv \frac{1}{2}\Big (
\rho_f\|\bu^i\|^2_{L^2(\Omega_f)}
+\rho_p \|w_t^i\|^2_{L^2(\omega_p)}+c_p\|p_p^i\|^2_{L^2(\Omega_p)}
+\|\bfeta^i\|^2_{E}+c_b\|p_b^i\|^2_{L^2(\Omega_b)}+D\|\Delta w^i\|^2_{L^2(\omega_p)}+{\gamma}\| w^i\|^2_{L^2(\omega_p)}
\Big )
$$
 and 
 $$k_b^i(t) \equiv k_b\big(c_bp^i(t)+\alpha_b\nabla \cdot \bfeta^i(t)\big).$$


We focus on the nonlinear term with the aim to obtain an estimate on $\overline{\mathbf s}$ through the functional $\mathcal E$. 
We can rewrite the nonlinear term as
\begin{align}
\int_0^T \left([k_b^1(t)\nabla p_b^1(t)-k_b^2(t)\nabla p_b^2(t)],\nabla \overline p_b\right)_{\Omega_b}dt = &~\int_0^T\big([k_b^1(t)-k_b^2(t)]\nabla p^1_b(t),\nabla \overline p_b(t)\big)_{L^2(\Omega_b)}dt \\
&+\int_0^T\big(k^2_b(t)\nabla \overline p_b(t),\nabla \overline p_b(t)\big)_{L^2(\Omega_b)}dt, \nonumber
\end{align}
where we used $\overline p = p^1_b-p^2_b$.
Using the lower bound on the permeability function $0<k_{\rm min}\le k_b$, and discarding the additional dissipation terms, we obtain the inequality
\begin{align}\nonumber
\mathcal E(\overline{\mathbf s}(T)) +&{{k_{\rm min}}} \int_0^T||\nabla \overline p_b||^2_{L^2(\Omega_b)}dt \le ~  \mathcal E(\overline{\mathbf s}(0))-\int_0^T\big([k_b^1(t)-k_b^2(t)]\nabla p^1_b(t),\nabla \overline p_b(t)\big)dt.
\end{align}

The nonlinear term on the right hand-side can be estimated in two ways, yielding the two cases in this theorem. 
In both cases we will use the Lipschitz hypothesis on the function $k_b$. 
Namely, let $L_k>0$ be the global Lipschitz constant so that
$$|k^1_b(t)-k^2_b(t)| = \big|k_b\big(c_bp_b^1(t)+\alpha_b\nabla \cdot \bfeta^1(t)\big)-k_b\big(c_bp_b^2(t)+\alpha_b\nabla \cdot \bfeta^2(t)\big)| \le L_k|c_b\overline p_b(t)+\alpha_b \nabla \cdot \overline {\bfeta}(t) |.$$
By the Cauchy-Schwarz inequality, we get
$$\int_0^T\big([k_b^1(t)-k_b^2(t)]\nabla p^1_b(t),\nabla \overline p_b(t)\big)dt  
\le L_k\int_0^T||\nabla p_b^1||_{L^{\infty}(\Omega_b)} ||c_b\overline p_b(t)+\alpha_b \nabla \cdot \overline {\bfeta}(t) ||_{L^2(\Omega_b)}
||\nabla \overline p_b||_{L^2(\Omega_b)}dt.$$
Now we estimate the right hand side by using regularity hypothesis, i.e. $\nabla p^1_b\in L^2(0,T^*;\mathbb{L}^{\infty}(\Omega_b))$. Namely, {{from Young's and triangle inequalities}} we have the following:
\begin{align}\int_0^T\big([k_b^1(t)-k_b^2(t)]\nabla p^1_b(t),\nabla \overline p_b(t)\big)dt 
 \le
 \int_0^T  L_k||\nabla p_b^1||_{L^{\infty}(\Omega_b)}  ||c_b\overline p_b+\alpha_b \nabla \cdot \overline {\bfeta}(t) ||_{L^2(\Omega_b)}||\nabla \overline p_b||_{L^2(\Omega_b)}dt 
 \nonumber 
 \\
 \nonumber
\le ~ \dfrac{L_k^2}{{{2\epsilon}}}||\nabla p^1_b||^2_{L^{\infty}(0,T;L^{\infty}(\Omega_b))}
\left\{
 \int_0^T||c_b \overline p_b||^2_{L^2(\Omega_b)}dt
+\int_0^T  ||\alpha_b\nabla \cdot \overline{\bfeta}(t)||^2_{L^2(\Omega_b)}dt 
\right\}
+ {{\frac{\epsilon}{2}}} \int_0^T||\nabla \overline p_b||^2_{L^2(\Omega_b)}dt.
\end{align}
Choosing $\epsilon={{k_{\rm min}}}$ to absorb the pressure term into the left hand-side, 
and by bounding the elasticity term from above by the full $||\cdot||_E$ norm, 
we have the following estimate:
\begin{align}\int_0^T\big([k_b^1(t)-k_b^2(t)]\nabla p^1_b(t),\nabla \overline p_b(t)\big)dt  
\le~ \int_0^T  L_k||\nabla p_b^1||_{L^{\infty}(\Omega_b)}  ||c_b\overline p_b+\alpha_b \nabla \cdot \overline {\bfeta}(t) ||_{L^2(\Omega_b)}||\nabla \overline p_b||_{L^2(\Omega_b)}dt 
\nonumber 
\\
\nonumber
\le \int_0^T\dfrac{L_k[c_b+\alpha_b^2]}{2k_{\rm min}}||\nabla {{p^1_b(t)}}||^2_{L^{\infty}(\Omega_b)}
\left\{ c_b||\overline p_b(t)||^2_{L^2(\Omega_b)}+||\overline{\bfeta}(t)||^2_E \right\}dt 
+ \frac{k_{\rm min}}{2} \int_0^T||\nabla \overline p_b||^2_{L^2(\Omega_b)}dt.
\end{align}
{{After discarding the dissipation term}}, we get another Gr\"onwall-type estimate:
\begin{equation}\label{Gron2}
\mathcal E(\overline{\mathbf s}(T)) \le C \mathcal E(\overline{\mathbf s}(0))+C\dfrac{L_k^2[c_b+\alpha_b^2]}{2k_{\rm min}}\int_0^T||\nabla p^1_b(t)||^2_{L^{\infty}(\Omega_b)}\mathcal E(\overline{\mathbf s}(t))dt 
\end{equation}
Uniqueness follows by applying the $L^2$-version of Gr\"onwall's inequality \cite{gronish}.

\end{proof}
We conclude this section by mentioning that in 3D, hypothesis of Proposition \ref{Uniq} would be satisfied for instance, if $p_b^1 \in L^2(0,T; H^{2.5+\epsilon}(\Omega_b))$.

\section{Acknowledgements} 
L. Bociu was partially supported by NSF-DMS 1555062 (CAREER). S. \v{C}ani\'{c} was partially supported by 
NSF-DMS 1853340, and by NSF-DMS 2011319.  B. Muha was partially supported by the Croatian Science Foundation project IP-2019-04-1140. J.T. Webster was partially supported by NSF-DMS 1907620.

B. Muha and J.T. Webster thank UMBC's Office of Research Development for the generous support through a SURFF grant, as well as the Department of Mathematics and Statistics for facilitating their collaboration in October and November, 2020.


\begin{thebibliography}{1}






\bibitem{AEN18}
Ilona Ambartsumyan, Vincent~J. Ervin, Truong Nguyen, and Ivan Yotov.
\newblock A nonlinear {S}tokes-{B}iot model for the interaction of a
  non-{N}ewtonian fluid with poroelastic media.
\newblock {\em ESAIM Math. Model. Numer. Anal.}, 53(6):1915--1955, 2019.

\bibitem{AKY18}
Ilona Ambartsumyan, Eldar Khattatov, Ivan Yotov, and Paolo Zunino.
\newblock A {L}agrange multiplier method for a {S}tokes-{B}iot
  fluid-poroelastic structure interaction model.
\newblock {\em Numer. Math.}, 140(2):513--553, 2018.

\bibitem{Auriault97}
Jean-Louis Auriault.
\newblock Poroelastic media.
\newblock In {\em Homogenization and porous media}, volume~6 of {\em
  Interdiscip. Appl. Math.}, pages 163--182, 259--275. Springer, New York,
  1997.
  
\bibitem{reg1} Avalos, G., Geredeli, P.G. and Muha, B., 2020. Wellposedness, spectral analysis and asymptotic stability of a multilayered heat-wave-wave system. Journal of Differential Equations,269(9), pp.7129-7156.

\bibitem{badia2009coupling}
Santiago Badia, Annalisa Quaini, and Alfio Quarteroni.
\newblock Coupling biot and navier--stokes equations for modelling
  fluid--poroelastic media interaction.
\newblock {\em Journal of Computational Physics}, 228(21):7986--8014, 2009.

\bibitem{banks1}
\newblock H.T. Banks, K. Bekele-Maxwell, L. Bociu, M. Noorman, G. Guidoboni,
\newblock Sensitivity Analysis
in Poro-Elastic and Poro-Visco-Elastic Models with Respect to Boundary Data,
\newblock Quarterly of
Applied Mathematics 75 (2017) 697-735

\bibitem{banks2}
\newblock H.T. Banks, K. Bekele-Maxwell, L. Bociu, M. Noorman, G. Guidoboni,
\newblock Local sensitivity via
the complex-step derivative approximation for 1-D poro-elastic and poro-visco-elastic models
\newblock Mathematical Control and Related Fields 9(4) (2019) 623-642

\bibitem{barucq1}
\newblock H. Barucq, M. Madaune-Tort, P. Saint-Macary,
\newblock Theoretical Aspects of Wave Propagation for Biot's Consolidation Problem,
\newblock Monografas del Seminario Matemtico Garca de Galdeano 31, 449--458 (2004) 449

\bibitem{barucq2}
\newblock H. Barucq, M. Madaune-Tort, P. Saint-Macary,
\newblock On nonlinear Biot's consolidation models,
\newblock Nonlinear Analysis 63 (2005) e985 -- e995


\bibitem{biot} Biot, M.A., 1941. General theory of three-dimensional consolidation. Journal of applied physics, 12(2), pp.155-164.

\bibitem{Biot55}
M.~A. Biot.
\newblock Theory of elasticity and consolidation for a porous anisotropic
  solid.
\newblock {\em J. Appl. Phys.}, 26:182--185, 1955.

\bibitem{Biot63}
M.~A. Biot.
\newblock Theory of stability and consolidation of a porous medium under
  initial stress.
\newblock {\em J. Math. Mech.}, 12:521--541, 1963.

\bibitem{BCLT} 
\newblock L. Bociu, L. Castle, I. Lasiecka, and A. Tuffaha, 
\newblock Minimizing Drag in A Moving Boundary Fluid-Elasticity Interaction, 
\newblock {\em Nonlinear Analysis} 197 (2020).

\bibitem{bgsw}
\newblock L. Bociu, G. Guidoboni, R. Sacco, and J. T. Webster,
\newblock Analysis of Nonlinear Poro-Elastic and Poro-Visco-Elastic Models,
\newblock {\em ARMA}, 222 (2016), pp. 1445--1519


\bibitem{MBE2}
\newblock L. Bociu, G. Guidoboni, R. Sacco,  and M. Verri,
\newblock On the role of compressibility in poroviscoelastic models,
\newblock Mathematical
Biosciences and Engineering, 16(5): 6167-6208, DOI: 10.3934/mbe.2019308

\bibitem{BW}
\newblock L. Bociu and J.T. Webster,
\newblock Nonlinear Quasi-static Poroelasticity
\newblock Accepted,{\em JDE}, 5/2021.

\bibitem{BN19}
\newblock L. Bociu and M. Noorman,
\newblock Poro-Visco-Elastic Models in Biomechanics: Sensitivity Analysis, 
\newblock Communications in Applied Analysis 23(1) (2019), 61-77

\bibitem{bukavc2012fluid}
\newblock M.~Buka{\v{c}}, I.~Yotov, and P.~Zunino.
\newblock An operator splitting approach for the interaction between a fluid
  and a multilayered poroelastic structure.
\newblock {\em Numerical Methods for Partial Differential Equations},
  31(4):1054--1100, 2015.
  
\bibitem{burman2014unfitted}
E. Burman and M. A. Fern{\'a}ndez,
\newblock An unfitted Nitsche method for incompressible fluid--structure interaction using overlapping meshes.
\newblock {\em Computer Methods in Applied Mechanics and Engineering}, 279 (2014), 497--514

 
\bibitem{guidoboni2006}
\newblock S. Canic, J. Tambaca, G. Guidoboni, A. Mikelic, C. J. 
Hartley, and D. Rosenstrauch
\newblock Modeling Viscoelastic Behavior of Arterial Walls and Their Interaction with Pulsatile Blood Flow,
\newblock {\em SIAM J. Appl. Math.}, 
67(1) (2006), pp. 164--193.


\bibitem{DESpaper}
\newblock S. Canic, Y. Wang, and M. Bukac
\newblock A next-generation mathematical model for drug-eluting stents
\newblock Submitted 2020.

\bibitem{cao}
\newblock Y. Cao, S. Chen, and A.J. Meir,
\newblock Analysis and numerical approximations of equations of nonlinear poroelasticity,
\newblock {\em DCDS-B}, 18 (2013), pp. 1253--1273.


\bibitem{causin2014poroelastic}
Paola Causin, Giovanna Guidoboni, Alon Harris, Daniele Prada, Riccardo Sacco,
  and Samuele Terragni.
\newblock A poroelastic model for the perfusion of the lamina cribrosa in the
  optic nerve head.
\newblock {\em Mathematical biosciences}, 257:33--41, 2014.

\bibitem{Ce17}
Aycil Cesmelioglu.
\newblock Analysis of the coupled {N}avier-{S}tokes/{B}iot problem.
\newblock {\em J. Math. Anal. Appl.}, 456(2):970--991, 2017.



\bibitem{coussy} Coussy, O., 2004. Poromechanics. John Wiley \& Sons.

\bibitem{cowin1999bone}
Stephen~C Cowin.
\newblock Bone poroelasticity.
\newblock {\em Journal of biomechanics}, 32(3):217--238, 1999.

\bibitem{detournay1993fundamentals}
Emmanuel Detournay and Alexander H-D Cheng.
\newblock Fundamentals of poroelasticity.
\newblock In {\em Analysis and design methods}, pages 113--171. Elsevier, 1993.

\bibitem{gronish} Dragomir, S.S., 2003. Some Gronwall type inequalities and applications. Nova Science.

\bibitem{dreher} Dreher, M. and J\" ungel, A., 2012. Compact families of piecewise constant functions in Lp (0, T; B). Nonlinear Analysis: Theory, Methods \& Applications, 75(6), pp.3072-3077.

\bibitem{FM03}
J.~L. Ferr\'{\i}n and A.~Mikeli\'{c}.
\newblock Homogenizing the acoustic properties of a porous matrix containing an
  incompressible inviscid fluid.
\newblock {\em Math. Methods Appl. Sci.}, 26(10):831--859, 2003.

\bibitem{Fung}
Y. C. Fung.
\newblock{\em Foundation of Solid Mechanics}.
\newblock Prentice-Hall, 1965.

\bibitem{huyghe1992porous}
Jacques~M Huyghe, Theo Arts, Dick~H Van~Campen, and Robert~S Reneman.
\newblock Porous medium finite element model of the beating left ventricle.
\newblock {\em American Journal of Physiology-Heart and Circulatory
  Physiology}, 262(4).
  
  \bibitem{Hsu}
  C.T. Hsu and P. Cheng.
  \newblock Thermal dispersion in a porous media.
  \newblock Int. J. Heat Mass Tran., 33(8)(1990), 1587--1597.

\bibitem{JagMik96}
Willi J\"{a}ger and Andro Mikeli\'{c}.
\newblock On the boundary conditions at the contact interface between a porous
  medium and a free fluid.
\newblock {\em Ann. Scuola Norm. Sup. Pisa Cl. Sci. (4)}, 23(3):403--465, 1996.

\bibitem{JagMik00}
Willi J\"{a}ger and Andro Mikeli\'{c}.
\newblock On the interface boundary condition of {B}eavers, {J}oseph, and
  {S}affman.
\newblock {\em SIAM J. Appl. Math.}, 60(4):1111--1127, 2000.

\bibitem{CzoMik15}
Anna Marciniak-Czochra and Andro Mikeli\'{c}.
\newblock A rigorous derivation of the equations for the clamped
  {B}iot-{K}irchhoff-{L}ove poroelastic plate.
\newblock {\em Arch. Ration. Mech. Anal.}, 215(3):1035--1062, 2015.

\bibitem{mikelic} 
A.  Marciniak-Czochra, Andro Mikelic.
A rigorous derivation of the equations for the clamped Biot-Kirchhoff-Love Poroelastic Plate.
\newblock{\em Arch. Rat. Mech. Anal.}, 215: 1035--1062, 2014.

\bibitem{MikTam16}
Andro Mikeli\'{c} and Josip Tamba\v{c}a.
\newblock Derivation of a poroelastic flexural shell model.
\newblock {\em Multiscale Model. Simul.}, 14(1):364--397, 2016.

\bibitem{MikWhe12}
Andro Mikeli\'{c} and Mary~F. Wheeler.
\newblock On the interface law between a deformable porous medium containing a
  viscous fluid and an elastic body.
\newblock {\em Math. Models Methods Appl. Sci.}, 22(11):1250031, 32, 2012.

\bibitem{reg2}   Muha, B., 2015. A note on optimal regularity and regularizing effects of point mass coupling for a heat-wave system. Journal of Mathematical Analysis and Applications, 425(2), pp.1134-1147.

\bibitem{ARMA} Muha, B. and Canic, S., 2013. Existence of a weak solution to a nonlinear fluid-structure interaction problem modeling the flow of an incompressible, viscous fluid in a cylinder with deformable walls. Archive for rational mechanics and analysis, 207(3), pp.919-968.

\bibitem{BorisSunNWE}
Muha, B. and Canic, S., 2014. Existence of a solution to a fluid-multi-layered-structure interaction problem. Journal of Differential Equations, 256(2), pp.658-706.

\bibitem{nitsche1971variationsprinzip}
J. Nitsche,
\newblock {\"U}ber ein Variationsprinzip zur L{\"o}sung von Dirichlet-Problemen bei Verwendung von Teilr{\"a}umen, die keinen Randbedingungen unterworfen sind,
\newblock {\em Abhandlungen aus dem mathematischen Seminar der Universit{\"a}t Hamburg}, 36 (1), 9--15, 1971

\bibitem{owc}
Sebastian Owczarek, 
\newblock A Galerkin method for Biot consolidation model. 
\newblock {\em Mathematics and mechanics of solids}, 15(1), pp.42-56, 2010.

\bibitem{RR}
\newblock{M. Renardy and R. C. Rogers}
\newblock An Introduction to Partial Differential Equations, 
\newblock Springer-Velag, 1992.

\bibitem{rohan1} 
\newblock E. Rohan and V. Luke\v{s},
\newblock Homogenization of the vibro-acoustic transmission on perforated plates,
\newblock{\em Appl. Math. Comput.}, (361), pp.821--845, 2019.

\bibitem{rohan2}
\newblock E. Rohan, and S. Naili,
\newblock Homogenization of the fluid-structure interaction in acoustics of porous media perfused by viscous fluid. 
\newblock {\em Z. Angew. Math. Phys.} 71, 137 (2020).

\bibitem{showmono}
\newblock R.E. Showalter,
\newblock Monotone Operators in Banach Space and Nonlinear Partial Differential Equations,
\newblock {\em AMS}, Mathematical Surveys and Monographs, 49, 1996. 

\bibitem{Sho05}
Ralph~E. Showalter.
\newblock Poroelastic filtration coupled to {S}tokes flow.
\newblock In {\em Control theory of partial differential equations}, volume 242
  of {\em Lect. Notes Pure Appl. Math.}, pages 229--241. Chapman \& Hall/CRC,
  Boca Raton, FL, 2005.
  
  \bibitem{show1} 
\newblock R.E. Showalter,
\newblock Diffusion in poroelastic media,
\newblock {\em JMAA}, 251 (2000), pp. 310--340.

\bibitem{DesaiPancreasReview}
S. Song, G. Faleo, R. Yeung, R. Kant, A.M. Posselt, T.A. Desai, Q. Tang, and S. Roy. Silicon nanopore membrane (SNM) for islet encapsulation and immunoisolation under convective transport. {\em Scientific reports}, 6(1), pp.1--9, 2016.

\bibitem{show2} 
\newblock N. Su and R.E. Showalter,
\newblock Partially saturated flow in a poroelastic medium,
\newblock {\em DCDS-B}, 1 (2001), pp. 403--420.

\bibitem{Tr}
\newblock F. Troltzsch, 
\newblock Optimal Control of Partial Differential Equations - Theory, Methods and Applications,
\newblock  AMS, 2010.


\bibitem{MBE} 
\newblock M. Verri, G. Guidoboni, L. Bociu, and R. Sacco,
\newblock The role of structural viscoelasticity in deformable
porous media with incompressible constituents: Applications in biomechanics. 
\newblock Mathematical
Biosciences and Engineering, 15(1551-0018 2018 4 933):933, 2018.

\bibitem{zenisek}
\newblock A. Zenisek,
\newblock The existence and uniqueness theorem in Biot's consolidation theory, 
\newblock {\em Appl. Math.}, 29 (1984), pp. 194--211.








\end{thebibliography}
\end{document}